\documentclass [12pt]{article}
\usepackage {amsfonts}
\usepackage {mathrsfs}
\usepackage {amsmath}
\usepackage {latexsym}
\usepackage {amssymb}
\usepackage {amsthm}
\usepackage {amscd}
\date{}
\title{\bf Finiteness of  Nichols algebras and Nichols (braided) Lie algebras}
\author{ Weicai Wu $^{a}$, Shouchuan Zhang $^{b}$  and Yao-Zhong Zhang $^{c}$   \\
\small\it $a$. School of Mathematics,  Hunan Institute of Science and Technology, Yueyang 414006,   China \\
\small \it $b$. Department  of Mathematics,   Hunan University, Changsha  410082,   China \\
\small \it $c$. School of Mathematics and Physics,   The University of Queensland, Brisbane 4072,   Australia\\
\small {\tt weicaiwu@hnu.edu.cn (WW); z9491@sina.cn (SZ); yzz@maths.uq.edu.au (YZZ)}
 }

\begin{document}
\newtheorem{Proposition}{Proposition}[section]
\newtheorem{Theorem}[Proposition]{Theorem}
\newtheorem{Definition}[Proposition]{Definition}
\newtheorem{Corollary}[Proposition]{Corollary}
\newtheorem{Lemma}[Proposition]{Lemma}
\newtheorem{Example}[Proposition]{Example}
\newtheorem{Remark}[Proposition]{Remark}

\maketitle 

\begin {abstract}  It is shown that if $\mathfrak B(V) $ is connected Nichols algebra of diagonal type with $\dim V>1$, then $\dim (\mathfrak L^-(V)) = \infty$
$($resp. $ \dim (\mathfrak L(V)) = \infty $$)$ $($ resp. $ \dim (\mathfrak B(V)) = \infty $$)$  if and only if $\Delta(\mathfrak B(V)) $ is an arithmetic root system and the quantum numbers (i.e. the fixed parameters) of generalized Dynkin diagrams of $V$ are of finite order. Sufficient and necessary conditions for $m$-fold adjoint action in $\mathfrak B(V)$ equal to zero, viz. $\overline{l}_{x_{i}}^{m}[x_{j}]^ -=0$ for $x_i,~x_j\in \mathfrak B(V)$,  are given.


\end {abstract}

\section {Introduction}\label {s0}
The question of finite-dimensionality of Nichols algebras dominates an important part of the recent developments in the theory of (pointed) Hopf algebras(see e.g. \cite {AHS08, AS10, An11, He05, He06a, He06b, WZZ15a, WZZ15b}. The interest in this problem comes from the lifting method of Andruskiewitsch and Schneider in classifying finite dimensional (Gelfand-Kirillov) pointed Hopf algebras  which are generalizations of quantized enveloping algebras of semi-simple Lie algebras.

Classification of arithmetic root systems is obtained in \cite {He05} and \cite {He06a}.
It is shown in \cite {WZZ15b} that   Nichols algebra $\mathfrak B(V)$  is
finite-dimensional if and only if Nichols braided Lie algebra $\mathfrak L(V)$ is finite-dimensional.

The main motivation for the subject of this paper and papers \cite {He05, He06a} is to answer the {\it Question 5.9} in \cite {A02} by Andruskiewitsch:
{\it Given a braided vector space $V$ of diagonal type, determine
when $\mathfrak B(V)$ is finite dimensional. If so, compute $\dim \mathfrak B(V )$, and give a nice presentation in terms of generators and relations.}

The first part of this problem was solved in papers \cite {He05} and \cite {He06a} under the following additional property (see \cite [Section 2.2] {He05}):
\vskip.1in
{\it
\noindent (P) The height of a PBW generator (i.e a hard super-letter) of $\mathbb Z^n$-degree $e$ is finite
if and only if $2 \le {\rm ord } (\chi (e, e)) < \infty $,  and in this case it coincides with ${\rm ord } (\chi (e, e))$.}
\vskip.1in

In this paper we show that the property (P) holds in any arithmetic root system (see Theorem \ref {2.4} below). Consequently, this paper solves  the first part of the above mentioned problem 5.9 for arithmetic root systems.

We will focus on when a Nichols algebra and a Nichols (braided) Lie algebra have finite dimensions. Such finite-dimensional Nichols algebras play a fundamental role in various subjects such as pointed Hopf algebras and logarithmic quantum field theories.   There exists a large number of  examples of infinite dimensional Lie algebras which are Nichols Lie algebras.

The main results of this paper are as follows: {\rm (i)} We show that if $\mathfrak B(V) $ is a connected Nichols algebra of diagonal type with $\dim V>1$, then $ \dim (\mathfrak L^-(V)) = \infty $
$($resp. $ \dim (\mathfrak L(V)) = \infty $$)$ $($ resp. $ \dim (\mathfrak B(V)) = \infty $$)$  if and only if $\Delta(\mathfrak B(V)) $ is an arithmetic root system and the quantum numbers (i.e. the fixed parameters  ) of the generalized Dynkin diagrams of $V$ are of finite order.
{\rm (ii)} We give the sufficient and necessary conditions for $m$-fold adjoint action in $\mathfrak B(V)$ equal to zero, viz. $\overline{l}_{x_{1}}^{m}[x_{2}]^ -=0$ and $\overline{l}_{x_{2}}^{m}[x_{1}]^ -=0$.

\section*{ Preliminaries}\label {s1}

Throughout this paper  braided vector space $V$ is of diagonal type with basis
$x_1,  x_2,  \cdots,  x_n$ and braiding $C(x_i \otimes x_j) = p_{ij} x_j \otimes x_i$ without special announcement.  Nichols algebra $\mathfrak B(V)$ over braided  vector space $V$¡¡ has two gradations.
A gradation is  the length such that $\mid u \mid = r$ of length of $u$ for $u = x_{i_1} x_{i_2} \cdots x_{i_r}$. Another gradation is $\mathbb Z^n$-graded such that ${\rm deg } (x_i) = e_i$ for $1\le i \le n$, where $E = \{ e_1, e_2, \cdots, e_n\}$ is a basis of $\mathbb Z^n.$

Define  linear map $p$ from $\mathfrak B(V) \otimes \mathfrak B(V) $ to $F$ such that
$p(u\otimes v) = \chi (deg (u),   deg (v)),  $ for any homogeneous element $u,   v \in \mathfrak B(V).$
For convenience,   $p(u\otimes v)$ is denoted by $p_{uv}$.  Let $\widetilde{p}_{uv} := p_{uv}p_{vu}.$

Let $A =: \{x_1,  x_2,  \cdots,  x_n\}$ a alphabet; $A^*$ denotes the set of all of  words in $A$; $A^+ =: A^*\setminus 1$.  Define $x_1 < x_2< \cdots < x_n$ and the order on $A^*$ is the lexicographic orderings. The concepts of words refer \cite {Lo83}.

\begin {Definition} \label {0.1} \cite [Def. 1] {Kh99}
A word $u$ is called a Lyndon word if $\mid u\mid =1$ or $\mid u\mid \geq2$, and for each
representation $u=u_1u_2$, where $u_1$and $u_2$ are nonempty
words, the inequality $u<u_2u_1$ holds.
\end {Definition}

 Any word $u\in A^*$ has a unique decomposition into the product of a non-increasing sequence of Lyndon words by \cite [Th.5.1.5] {Lo83}, If $u$ is a Lyndon word with $\mid u \mid >1$, then there uniquely exist two Lyndon words $v$ and $w$ such that $u =vw$ and $v$ is shortest (see \cite [Pro. 5.1.3]{Lo83} and \cite {He07}) (the composition is called the Shirshov decomposition of $u$).

\begin {Definition} \label {0.2} We inductively define a linear map $[\  \ ]$ from $A^+$ to $\mathfrak B(V)$  as follows:

(1) $[u] =: u$  when  $u$ is a letter;

(2) $[u] =: [w][v]-p_{wv}[v][w]$ when $u$ is a Lyndon word with $\mid u \mid >1$ and $u =vw$ is a Shirshov decomposition;

(3) $[u] =: [[[l_1, l_2], l_3] \cdots l_n] $, when $u = l_1l_2 \cdots l_n$ is a Shirshov decomposition.
\end {Definition}

$[u]$ is called a nonassociative word for any $u \in A^+$. $[u]$ is called   a super-letter if $u $ is a  Lyndon word. A multiplication of super-letters is called a super-word.

\begin {Definition} \label {0.3} \cite [Def. 6] {Kh99}
A super-letter $[u]$ is said to be hard if it is not a linear
combination of super-words with length $\mid u \mid $ in greater super-letters than $[u]$.
\end {Definition}
\begin {Definition} \label {1.4} \cite [Def 7]{Kh99} or \cite [before Th. 10] {He07} We said that the height of a super-letter $[u]$ with length  $d$ equals a natural number $h$ if $h$ is least with the following properties:

(1) $p_{uu}$ is a $t$-primitive root of unity,   and $h=t$;

(2) super-word $[u]^h$ is a linear combination of super-words of length $hd$ in greater super-letters than $[u].$

If the number $h$ with above properties does not exist then we say that the height of $[u]$ is infinite.

\end {Definition}

Let ${\rm ord} (p_{uu})$ denote the  order of $p_{uu}$ with respect
to multiplication.

$D(V)=D =: \{[u] \mid [u] \hbox { is a hard super-letter }\}$,we know $h_{u}>1$ for $\ \forall \ [u]\in D$.

$P(V)= P = :\{[u_{1}]^{k_1}[u_{2}]^{k_2}\cdots  [u_{s}]^{k_s}\  |\  [u_{i}]\in D,  k_i,s \in \mathbb N_0; 0 \le k_i < h_{u_i};  1\le i \le s;   u_s<u_{s-1}<\cdots< u_1\}$

If $[u] \in D$ and
${\rm ord }  (p_{u, u}) =m>1$ with $h_u = \infty$, then $[u]$ is called an $m$-infinity element.

\begin {Theorem} \label {1.5} \cite [Th. 2] {Kh99} or \cite [Th. 1.4.6] {He05},
$P$ is a basis of $\mathfrak B (V)$.
\end {Theorem}

Let $G$ be a nonempty set; $R$   a nonempty subset of  $G\times G$, and  $\circ: R \rightarrow G$ a
map of sets. The pair $(G, \circ )$ is called a groupoid if it satisfies the following
conditions.

${\rm GD_1}$  If $(x, y) \in R$ then each of the three elements $x, y, x \circ  y $ is uniquely
determined by the other two.

${\rm GD_2}$ If $(x, y),(y, z) \in R $ then $(x\circ y, z),(x, y\circ z) \in R $ and $(x\circ y)\circ z = x\circ (y\circ z).$

${\rm GD_3}$ If $(x, y),(x\circ y, z) \in R$ then $(y, z),(x, y\circ z) \in R$ and $(x\circ y)\circ z = x\circ (y\circ z)$.

${\rm GD_4}$ If $(y, z),(x, y\circ z) \in R$ then $(x, y),(x\circ y, z) \in R$ and $(x\circ y)\circ z = x\circ (y\circ z)$.

${\rm GD_5}$ For all  $ x \in G$ there exist unique elements $e, f, y \in G$ such that $(e, x),
(x, f),(y, x) \in R$, $e \circ  x = x \circ  f = x,$ and $y \circ  x = f$.

${\rm GD_6}$ If $e \circ  e = e, f \circ  f = f $ for certain $e, f \in G$ then there exists $x \in G$ such
that $e \circ  x = x \circ  f = x.$

Let $E$ and  $F$ be  bases of $\mathbb Z^n$
 and $\chi $ a bicharacter on $\mathbb Z^n.$
Define
$m (e', e'') :=  {\rm min} \{ m \in \mathbb N_0 \mid \chi (e', e') ^m \chi (e', e'')
\chi (e'', e') = 1 \ {\rm or } \ \chi (e', e') ^{m+1} =1,$ $ \chi (e', e') \not= 1 \}$ for $e', e''\in F$ with $e' \not= e''$ and $m (e', e') := -2$.
$s_{e', F} (e'') := e'' + m (e', e'') e' $ for $e', e'' \in F.$

Let $\widetilde{W}$ denote the groupoid consisting of all pairs $(T, E)$,  where $T \in
{\rm Aut}(\mathbb Z^n )$
 and $ E$ is a basis of $\mathbb Z^n$,  and the composition $(T_1, E_1) \circ  (T_2, E_2)$ is
defined (and is then equal to $(T_1T_2, E_2))$ if and only if $T_2(E_2) = E_1$ (i.e. $R = \{  (T_1, E_1) \times   (T_2, E_2) \mid  T_2(E_2) = E_1 \}$).

Let $E$ be a basis of $\mathbb Z^n$  and $\chi $ a bicharacter on $\mathbb Z^n.$
Define  $ W_{\chi, E}$ as the smallest subgroupoid of $\widetilde{W} $ which contains $({\rm id} , E)$, and if
$( {\rm id}, F) \in W_{\chi, E} $ for a basis $F$ of $\mathbb Z^n$
and $f \in F$, then $(s_{f,F} , F),({\rm id}, s_{f,F} (F)) \in
W_{\chi,E}$ whenever $s_{f,F}$ is defined. It is called the Weyl groupoid associated to the pair $(\chi, E)$.

The groupoid $W_{\chi,E}$ is called full,
if $s_{f,F}$ is well-defined for all bases $F$ of $\mathbb Z^n$
 with $( {\rm id}, F) \in W_{\chi,E}$ and for all $f \in F$.

A triple $(\Delta , \chi, E)$ is called an arithmetic root system
if $ W_{\chi,E} $ is full and finite wirh $\Delta   :=
\cup \{F \mid ({\rm id}, F) \in W _{\chi,E}\}.$

$\Delta ^+(\mathfrak B(V)): =  \{ \deg (u) \mid [u]\in D\}$ and
$ \Delta (\mathfrak B(V)) := \Delta ^+(\mathfrak B(V)) \cup \Delta ^-(\mathfrak B(V))$. Let  $\mathfrak L(V)$ denote the  braided Lie algebras generated by $V$ in $\mathfrak B(V)$ under Lie operations $[x, y]=yx-p_{yx}xy$,  for any homogeneous elements $x, y\in \mathfrak B(V)$. $(\mathfrak L(V), [\ ])$ is called Nichols braided Lie algebra of $V$. Let  $\mathfrak L^-(V)$ denote the  Lie algebras generated by $V$ in $\mathfrak B(V)$ under Lie operations $[x, y]^-=yx-xy$,  for any homogeneous elements $x, y\in \mathfrak B(V)$. $(\mathfrak L^-(V), [\ ]^-)$ is called Nichols  Lie algebra of $V$.

\vskip.1in
The dual $\mathfrak B(V^*) $ of Nichols algebra $\mathfrak B(V) $ of rank $n$
in \cite [Section 1.3]{He05}. Let $y_{i}$ be a dual basis of $x_{i}$. $\delta (y_i)=g_i^{-1} \otimes y_i$,  $g_i\cdot y_j=p_{ij}^{-1}y_j$ and $\Delta (y_i)=g_i^{-1}\otimes y_i+y_i\otimes 1.$ There exists a bilinear map
$<\cdot, \cdot >$: $(\mathfrak B(V^*)\# FG)\times\mathfrak B(V)$ $\longrightarrow$ $\mathfrak B(V)$ such that

$<y_i, uv>=<y_i, u>v +g_i^{-1}.u<y_i, v>$ and $<y_i, <y_j, u>>=<y_iy_j, u>$

\noindent for any $u, v\in\mathfrak B(V)$. Furthermore,  for any $u\in \oplus _{i=1}^\infty \mathfrak B(V)_{(i)}$,
one has that  $u=0$ if and only if $<y_i, u>=0$ for any $1\leq i \leq n.$

We have the braided Jacobi identity as follows:
\begin {eqnarray} [[u, v], w]&=&[u, [v, w]]+p_{vw}^{ -1} [[u, w], v]+(p_{wv}-p_{vw}^{-1})v\cdot[u, w],\label {e0.1}\\
{} [u, v \cdot w]&=&p_{wu}[uv]\cdot w + v \cdot [uw].\label {e0.2}
\end {eqnarray}

\begin {Theorem} \label {1.6}(\cite [Th. 2.5.3] {He05}) If $\chi $ is a bicharacter on $\mathbb Z^n$ and $(\Delta, \chi , E)$ is an arithmetic
root system then for the braided vector space $V$ of diagonal type with
$\dim V = n$ and with structure constants $q_{ij} := \chi (e_i
, e_j )$ one has $\Delta(\mathfrak B(V)) =
\Delta$. Conversely, if $V$  is a braided vector space of diagonal type such that
$\Delta^+
(\mathfrak B(V)) $ is finite then $(\Delta (\mathfrak B(V)), \chi , E) $ is an arithmetic root system, where $E= \{e_1, e_2, \cdots, e_n\}$ and
$\chi (e_i, e_j) =q_{ij}$ for $1\le i, j\le n.$
\end {Theorem}

If $(\Delta , \chi, E)$ is an arithmetic root system and
 a hyperplane  $H  $ with $0 \in H$ (i.e. $\dim H = n-1$), then there exists  $E_H := \{   f_1, f_2 , \cdots, f_l\}$ such that
$(\Delta \cap H,  \chi _ {\Gamma \times \Gamma} , E_H)$ is  an arithmetic root subsystem
with $\Gamma = \mathbb R(\Delta \cap H) \cap \mathbb Z^n, $ written  $ \Delta ( \chi, f_1, f_2, \cdots, f_l) = (\Delta \cap H,  \chi _ {\Gamma \times \Gamma} , E_H)$ in short (see \cite [Pro. 2.7.1]{He05}).

Throughout,  $\mathbb Z =: \{x \mid  x \hbox { is an integer}\}.$ $\mathbb R =: \{ x \mid x \hbox { is a real number}\}$.
$\mathbb N_0 =: \{x \mid  x \in \mathbb Z, x\ge 0\}.$
$\mathbb N =: \{x \mid  x \in \mathbb Z,  x>0\}$. $\mathbb R _+ =: \{ x \ge 0 \mid x \hbox { is a real number}\}$. $F$ denotes the base field,   which is an algebraic closed field with  characteristic zero. $F^{*}=F\backslash\{0\}$. $\mathbb S_{n}$ denotes symmetric group, $n\in\mathbb N$.

\section {$m$-infinity elements }
In this section we show that there does not exist any $m$-infinity elements in arithmetic root systems.
\begin {Lemma} \label {2.1}  Assume that $\widetilde{f_1}, \widetilde{f_2}, \cdots, \widetilde{f_{\widetilde{l}}}$ are linearly independent  in  arithmetic root system $(\Delta, \chi, E )$ with $\widetilde{ H} := \mathbb R (  \widetilde{f_1}, \widetilde{f_2}, \widetilde{f_3}, \cdots, \widetilde{f_{\widetilde{l}}})$. If
\begin {eqnarray}\label {ep2.2.1} \Delta \cap \widetilde{H} \subseteq \mathbb R_+ ( \widetilde{f_1}, \widetilde{f_2}, \cdots, \widetilde{f_{\widetilde{l}}} ) \cup -\mathbb R_+ ( \widetilde{f_1}, \widetilde{f_2}, \cdots, \widetilde{f_{\widetilde{l}}} ), \end {eqnarray}
then there exists a hyperplane  $H  $ with $0 \in H$  and $E_H = \{   f_1, f_2 , \cdots, f_l\}$ with $l = \widetilde{l}$, such that  $\Delta \cap H = \Delta \cap \widetilde{H}$ and
$(\Delta \cap H,  \chi _ {\Gamma \times \Gamma} , E_H)$ is  an arithmetic root system
with $\Gamma = \mathbb R(\Delta \cap H) \cap \mathbb Z^n. $
\end {Lemma}
\noindent {\it Proof.} It is clear that there exists $H_1$ with $H_1 \cap \Delta = \emptyset$ such that
$H = \widetilde{H} \oplus H_1$ with $\dim H = n-1.$ By \cite [Pro. 2.7.1]{He05},
there exists  $E_H = \{   f_1, f_2 , \cdots, f_l\}$ such that
$(\Delta \cap H,  \chi _ {\Gamma \times \Gamma} , E_H)$ is  an arithmetic root system
with $\Gamma = \mathbb R(\Delta \cap H) \cap \mathbb Z^n $  and \begin {eqnarray}\label {ep2.2.2}E_H \subseteq \Delta \cap H =  \Delta \cap \widetilde{H}\subseteq \mathbb R E_H.\end {eqnarray} Considering (\ref {ep2.2.1}) and (\ref {ep2.2.2}), we have $l = \widetilde{l}$.\hfill $\Box$

\begin {Lemma} \label {2.2}  If  $(\Delta(\mathfrak B(V)), \chi, E) $ is an arithmetic root system with   $[v],  [w]\in D$ and ${\rm deg }(v) \not= {\rm deg }(w),  $  then there exists an  arithmetic root subsystem  $\Delta ( \chi; f_1,  f_2)$ of  $\Delta(\mathfrak B(V)) $ such that ${\rm deg } (v)$ and ${\rm deg } (w)$ are in the subsystem. Furthermore,  if $[u] = [[v],    [w]]$ is the Shirshov decomposition of $[u]\in D$  then the subsystem  $\Delta ( \chi; f_1,  f_2)$ is connected and  ${\rm deg } (u)$ are in the subsystem.
\end {Lemma}

\noindent {\it Proof.} We show this  by the following steps.

{\rm (i)} Let $v_1 = v$  and $v_2 = w$. If  $\deg(v)-k\deg(w) \notin  \mathbb N  \cdot \Delta(\mathfrak B(V))$ for any $k \in \mathbb N$,  let $\widetilde{f_1} := {\rm deg } (v_1)$  and $\widetilde{f_2} := {\rm deg } (v_2)$.
If  there exists a $k_1 \in \mathbb N$ such that $\deg(v_1)-k_1\deg(v_2)  \in \mathbb N \cdot \Delta(\mathfrak B(V))$,  then there exist  $ v_3 \in \Delta^+(\mathfrak B(V))$ and  $l_1 \in \mathbb Z\setminus \{0\}$ such that $\deg(v_1)-k_{1}\deg(v_2) = l_1 \deg(v_{3})$. Keeping on the step,  we obtain   $v_i \in \Delta(\mathfrak B(V))$,  $k_i \in \mathbb N$ and $ l _i \in  \mathbb Z \setminus \{0\}$ such that $\deg(v_i)-k_{i}\deg(v_{i+1}) =  l _{i}\deg(v_{i+2})\in \mathbb N\cdot\Delta (\mathfrak B(V))$ for $i =1,  2,  \cdots$.

 {\rm (ii)} $\deg(v_i) \not= \deg(v_j)$ when $i\not= j.$ Indeed,  for any $t$,  there exist $r_t,  s_t \in \mathbb R$ such that  $\deg(v_t) = r_t \deg(v_1) +s_t\deg( v_{2})$ by Part (i). If $\deg(v_i) = \deg ( v_j)$,  then $ (r_i -r_j)\deg(v_1) = \deg(v_2) (s_j -s_i)$,  which contradicts \cite [Cor. 2.5.4 ] {He05}.

{\rm (iii)} By Part {\rm (ii)},  the set  $\{ v_i \mid i=1,  2,  \cdots \}$ in Part {\rm (i)} is finite since $\Delta(\mathfrak B(V))$ is finite.

{\rm (iv)} By Part {\rm (i)}. {\rm (ii)} and {\rm (iii)},  there exists $i_0 \in \mathbb N$ such that  $\deg(v_i)-k_{i}\deg(v_{i+1}) =  l _{i}\deg(v_{i+2})\in \mathbb N\cdot\Delta (\mathfrak B(V))$ for $1 \le i < i_0$ and $\deg(v_{i_0})-k\deg(v_{i_0+1})\notin \mathbb N\cdot\Delta(\mathfrak B(V))$ for any $k \in \mathbb N$. Let $\widetilde{f_1} := {\rm deg } (v_{i_0})$  and $\widetilde{f_2} := {\rm deg } (v_{i_0+1})$. Let $\widetilde{H} :=\mathbb R \{\widetilde{f _1},  \widetilde{f_2}\} = {\rm span } \{\widetilde{f _1},  \widetilde{f_2}\}$;
 By \cite [  Lemma 2.7.2] {He05} and Lemma \ref {2.1}, there exists  an arithmetic root subsystem
$\Delta (\chi; f_1, f_2)$ such that ${\rm deg } (v)$ and ${\rm deg } (w)$ are in  $\Delta (\chi; f_1, f_2)$.

{\rm (v)}  The subsystem  $\Delta ( \chi; f_1,  f_2)$  is connected when $[u] = [[v],    [w]]$ is the Shirshov decomposition of $[u]\in D$. Otherwise,
  $ {\rm deg } (v),   {\rm deg } (w),      {\rm deg } (u)\in \{ f_1,  f_2 \}$,  which   is a contradiction.

\hfill $\Box$

\begin {Proposition} \label {2.3} If $(\Delta , \chi, E)$ is an arithmetic root system, then for every $\alpha \in \Delta\setminus E$, there exists a connected  arithmetic root subsystem
$(\Delta \cap H,  \chi _ {\Gamma \times \Gamma} , E_H)$ such that $\alpha \in \Delta \cap H$ and $\mid E_H \mid =2$
\end  {Proposition}
\noindent {\it Proof.}  It follows from Lemma \ref {2.2}.\hfill $\Box$

\begin {Theorem} \label {2.4} If $\mathfrak B(V) $ is connected Nichols algebra of diagonal type with $\dim V>1$ and $\Delta(\mathfrak B(V)) $ is an arithmetic root system,  then $p_{u, u}\neq1$ for $\forall\ [u]\in D$ and there do not exist any $m$-infinity elements.
\end {Theorem}

\noindent {\it Proof.} Set $[u] = [[v],    [w]]$ is the Shirshov decomposition of $[u]$. By Lemma \ref {2.2},  ${\rm deg} (u)$ is in  a connected subsystem with rank $2$. We complete the proof using \cite [A.2] {WZZ15b} and \cite [Th. 1.6.1] {He05}.  \hfill $\Box$

\section {Classification  of   Nichols  Lie algebras }

In this section we give the classification  of   Nichols  Lie algebras.

Let $l_{u}^0[v]:= [ v]$,  $l_{u}^i[v]:=  [[u], l_{u}^{i-1} [ v]].$ Similarly define $r_{u}^i[v].$
Let $\overline{l}_{u}^0[v]^-:= [ v]^-$,  $\overline{l}_{u}^i[v]^-:=  [[u]^-, \overline{l}_{u}^{i-1} [ v]^-]^-.$ Similarly define $\overline{r}_{u}^i[v]^-.$ In fact,
$\overline{l}_{u}^i[v]^-:= [u, [u, \cdots, [u, v]^-\cdots ]^-]^-$.

\begin {Lemma} \label {3.1} Assume that $\mathfrak B(V) $ is a  Nichols algebra of diagonal type. Let  $a:=p_{ii}^{-1}$,  $b:=p_{ij}^{-1}$,  $c:=p_{ji}^{-1}$ with $i\not= j$. Then

{\rm (i)} $\bar{l}_{i}^{m}[j]^-=0$,  when $b=1, c=1$.

{\rm (ii)} \begin {eqnarray}   <y_{i}^{k}, \bar{l}_{i}^{m}[j]^-> &=& \sum \limits _{l=0}^{k-1}(a^{l}-a^{m-1-l}b)<y_{i}^{k-1}, \bar{l}_{i}^{m-1}[j]^->\nonumber  \\
& &+a^{k}x_{i}<y_{i}^{k}, \bar{l}_{i}^{m-1}[j]^->
-<y_{i}^{k}, \bar{l}_{i}^{m-1}[j]^->x_{i} \ \ \hbox { for } k \le m, \label {e2.5}\\
<y_{i}^{m}, \bar{l}_{i}^{m}[j]^->&=&(1-b)^m(m)_{a}! x_j, \label {e2.5.2} \\
 <y_{j}y_{i}^{m}, \bar{l}_{i}^{m}[j]^->&=&(1-b)^m(m)_{a}!.   \label {e2.5.3}
\end {eqnarray}

{\rm (iii)} Assume that $a=1$. Then $\bar{l}_{i}^{m}[j]^-\not=0$,  when $b\not=1$ or $c\not=1$.

{\rm (iv)} Assume that $a\not=1$. Then $\bar{l}_{i}^{m}[j]^-\not=0$,  when $ {\rm ord }(a) > m$ with $c\not=1$ or $ b\not=1$.

{\rm (v)} $\dim (\mathfrak L^-(V))=\infty$ when there exists $i$ and $j$ with $i\not=j$,    $b\not=1$ or $c\not=1$ and ${\rm ord}(a)=1$ or $\infty$.
\end {Lemma}

\noindent {\it Proof.} {\rm (i)}  follows from  \cite [Lemma 4.1 {\rm (i)}(vi)] {WZZ15b}.

{\rm (ii)} We  show (\ref {e2.5})
 by induction on $k$ for $k \le m$. See \begin {eqnarray*} <y_{i}, \bar{l}_{i}^{m}[j]^-> &=& <y_{i},  x_i\bar{l}_{i}^{m-1}[j]^- -  \bar{l}_{i}^{m-1}[j]^-x_i > \\
&=& \bar{l}_{i}^{m-1}[j]^- + a x_i  <y_{i},  \bar{l}_{i}^{m-1}[j]^-> -  < y_i,  \bar{l}_{i}^{m-1}[j]^->x_i  - a^{m-1}b \bar{l}_{i}^{m-1}[j]^- \\
&= & (1- a ^{m-1}b ) \bar{l}_{i}^{m-1}[j]^- + a x_i < y_i,  \bar{l}_{i}^{m-1}[j]^-> - < y_i,  \bar{l}_{i}^{m-1}[j]^-> x_i.
 \end {eqnarray*}
Thus equation (\ref {e2.5}) holds when $k=1.$ Assume $k>1.$ See
 \begin {eqnarray*} <y_{i} ^k, \bar{l}_{i}^{m}[j]^->
 &=&
  <y_{i},  \sum \limits _{l=0}^{k-2}(a^{l}-a^{m-1-l}b)<y_{i}^{k-2}, \bar{l}_{i}^{m-1}[j]^-> \\
& &+ a^{k-1}x_{i}<y_{i}^{k-1}, \bar{l}_{i}^{m-1}[j]^->-<y_{i}^{k-1}, \bar{l}_{i}^{m-1}[j]^->x_{i} > \\
 & & ~~~~~~~~ \hbox { by inductive hypothesis}\\
&=& \sum \limits _{l=0}^{k-2}(a^{l}-a^{m-1-l}b)<y_{i}^{k-1}, \bar{l}_{i}^{m-1}[j]^-> \\
& &+ a^{k-1} <y_{i}^{k-1}, \bar{l}_{i}^{m-1}[j]^->  +a^{k-1} x_ia<y_{i}^{k}, \bar{l}_{i}^{m-1}[j]^->\\
& &-<y_{i}^{k}, \bar{l}_{i}^{m-1}[j]^->x_{i} - a ^{m-k} b<y_{i}^{k-1}, \bar{l}_{i}^{m-1}[j]^-> \\
&=& \hbox  { the right hand side of (\ref {e2.5})}.
\end {eqnarray*} Consequently,  (\ref {e2.5} ) holds.

Now we show (\ref {e2.5.2}) by induction on $m$.  (\ref {e2.5.2}) follows from  \cite [Lemma 4.1{\rm (i)}]{WZZ15b} when $m=1.$

By (\ref {e2.5}),  one obtains
\begin {eqnarray*}  <y_{i}^{m}, \bar{l}_{i}^{m}[j]^-> &=& \sum \limits _{l=0}^{m-1}(a^{l}-a^{m-1-l}b)<y_{i}^{m-1}, \bar{l}_{i}^{m-1}[j]^->\\
&=& (m)_a (1-b)<y_{i}^{m-1}, \bar{l}_{i}^{m-1}[j]^->\\
&=&   \hbox {the right hand side of  (\ref {e2.5.2})}.  \end {eqnarray*}
Therefore,  ( {\ref {e2.5.2}} ) and (\ref {e2.5.3}) hold.

{\rm (iii)} If $a=1$,  then  $<y_{i}^{m}y_{j}, \bar{l}_{i}^{m}[j]^->=(c-1)^m(m)!$ and $<y_{j}y_{i}^{m}, \bar{l}_{i}^{m}[j]^->=(1-b)^m(m)!$ by \cite [Lemma 4.1(v)] {WZZ15b} and {\rm (ii)}.

{\rm (iv)}  It follows from  {\rm (ii)} and  \cite [Lemma 4.1 {\rm (v)}] {WZZ15b}.

{\rm (v)} If $a=1$,  it is clear by {\rm (iii)}. If $a\not=1$,  then {\rm (v)} follows from  \cite [Lemma 4.1{\rm (iii)}] {WZZ15b} and {\rm (iv)}. \hfill $\Box$

\vskip.2in
Let $I$ be a graded ideal of $T(V)$ as algebras such that $\mathfrak B(V) = T(V)/I$.
$u$ is called a standard   word with respect to $\mathfrak B(V)$ if $u $ can not be  written as  linear combination of  strictly greater  words in $\mathfrak B(V)$.

Let $ {\rm SW} (\mathfrak B(V))=: \{ u \in A^* \mid  u \hbox { is a standard word with respect to } \mathfrak B(V) \}$, written as {\rm SW} in short;
 ${\rm L }:= \{ u\in A^* \mid u \hbox { is a Lyndon word}\}. $ Notice that we view
$ {\rm SW} (\mathfrak B(V))$ and ${\rm L }$ are in $\mathfrak B(V)$ often for convenience.
Let $ {\rm HBLLW}(\mathfrak B(V)) =: \{ u  \in {\rm L}\mid  [u] \hbox { is a hard super-letter }\}$.

\begin {Lemma}\label {3.2}
{\rm (i)} $ {\rm SW} $ is a basis of $\mathfrak B(V).$

{\rm (ii)} Any factor of a standard  word is a standard  word.

{\rm (iii)} If $u$ is a standard word, then  $u = u_1u_2 \cdots u_r$ with $u_! \ge u_2\ge \cdots \ge u_r $ and $u_i \in {\rm SW } \cap {\rm L}$ for $1\le i \le r.$

\end {Lemma}
\noindent {\it Proof.}  {\rm (i)} If a word $u \notin {\rm SW},$ then $u = \sum\limits _{v \in {\rm SW},  \mid v\mid = \mid u \mid, v >u} a_v v$. It is clear that ${\rm SW}$ is linearly independent in $\mathfrak B(V).$ Consequently, ${\rm SW}$ is a basis of $\mathfrak B(V).$

{\rm (ii)} It is clear.

{\rm (iii)} It follows from Part {\rm (ii)} and \cite [Th. 5.1.5] {Lo83}. \hfill $\Box$

\begin {Lemma}\label {5.2}
{\rm (i)}  If  $l \in {\rm L}, $  then   $ [l]^- =  l + \sum \limits _ {w > l, \mid l \mid = \mid w\mid} a_w w$ in $\mathfrak B(V),$  where $ a_w \in k$.

{\rm (ii)} If $l \in  {\rm L}, $ then  $ [l] =  a_l l + \sum \limits _ {w > l, \mid l \mid = \mid w\mid} a_w w$ in $\mathfrak B(V),$ where $a_l, a_w \in k$ with $a_l\not=0.$

{\rm (iii)}   $ {\rm SW} \cap  {\rm L} \subseteq  {\rm HBLLW}.$

\end {Lemma}
\noindent {\it Proof.}  {\rm (ii)} We show this by induction on $\mid l \mid.$  It is clear when $\mid l \mid =1$ since $[l] = l.$ Assume that $l = uv$ is the Shirshow decomposition of $l$.
 If $u' > u$ and  $v'>v$ with $\mid u'\mid = \mid u\mid $
  and $\mid v'\mid = \mid v\mid $, then $u'v' > uv =l$ and $v'u' >vu >l.$
 \begin {eqnarray*} [l] &=& [v] [u] - p _{v, u}[u][v]\\
 &=& (a_v' v + \sum \limits _ {v' > v, \mid v' \mid = \mid v \mid} a_{v'}' v')(a_u 'u + \sum \limits _ {u' > u, \mid u' \mid = \mid u \mid} a_{u'}' u') \\
&& - p _{vu}(a_u 'u + \sum \limits _ {u' > u, \mid u' \mid = \mid u \mid} a_{u'}' u')
 (a_v' v + \sum \limits _ {v' > v, \mid v' \mid = \mid v \mid} a_{v'}' v')\ \ ( \hbox {by inductive  hypothesis}) \\
&=&  a_ll + \sum \limits _ {w > l, \mid l \mid = \mid w\mid} a_ww.
\end {eqnarray*}

{\rm (i)}  The proof is similar to the proof of {\rm (ii)}.

{\rm (iii)} If $l  \notin {\rm  HBLLW}$ with $l \in {\rm L}$, then
 \begin {eqnarray*}  [l] &=& \sum \limits _{r=1} ^m  \sum \limits _ {u_i>l, 1\le  i \le r, \mid u_1u_2\cdots u_r \mid =  \mid l \mid  } a_{u}' [u_1] [u_2] \cdots [u_r] \\
&=&  \sum \limits _{r=1} ^m \sum \limits _ {u_i>l,  1\le  i \le r, \mid u_1u_2\cdots u_r \mid =  \mid l \mid } a_{u}' \prod _{i=1}^r ( a_{u_i}''u_i   + \sum \limits_{w_{ij} > u_i, \mid w_{ij} \mid =  \mid u_i \mid} a'' _{w_{ij}}  w_{ij}) \\
&=&\sum \limits _ {u>l, \mid u \mid =  \mid l \mid}   a_u u
\end {eqnarray*}
and $ [l]  =  a_l''' l + \sum \limits _ {w > l, \mid l \mid = \mid w\mid} a_w '''w$ by (ii).
Consequently, $ l =  \sum \limits _ {w > l, \mid l \mid = \mid w\mid} a_w w$ and  $l \notin  {\rm SW} \cap  {\rm L}.$.

\hfill $\Box$

\begin {Theorem}\label {5.3}   If  $\mathfrak B(V)$ is  a graded  braided Hopf algebra with diagonal braiding, then
${\rm SW } (\mathfrak B(V))\cap {\rm L} = {\rm HBLLW} (\mathfrak B(V)).$
Furthermore, if $\Delta ( \mathfrak B(V)  )$ is not  an arithmetic root system, then $\dim \mathfrak L^-( \mathfrak B(V)) = \infty.$

\end {Theorem}
\noindent {\it Proof.}  ${\rm SW }\cap {\rm L} \subseteq {\rm HBLLW}$ by Lemma \ref {5.2} {\rm (iii)}. If there exists $u \in  {\rm HBLLW}$ and $u\notin {\rm SW}\cap {\rm L}, $  then $u = \sum \limits_{v_1 \ge v_2 \ge \cdots \ge v_r, v_i \in {\rm SW }\cap {\rm L},
\mid v_1v_2\cdots v_r \mid = \mid u \mid, v_1v_2\cdots v_r>u} a_v v_1v_2 \cdots v_r$ by Lemma \ref {3.2} {\rm (iii)}. Consequently, $[u] =  \sum \limits_{v_1 \ge v_2 \ge \cdots \ge v_r, v_i \in {\rm SW }\cap {\rm L},
\mid v_1v_2\cdots v_r \mid = \mid u \mid, v_1v_2\cdots v_r>u} a_v [v_1][v_2] \cdots [v_r], $
 which contradicts to $[u], [v_1][v_2] \cdots [v_r] \in P.$
\hfill $\Box$

\begin {Proposition} \label {3.3} Assume that $\mathfrak B(V) $ is a connected  Nichols algebra of diagonal type. Then $\dim(\mathfrak B(V))=\infty$ and  $ \dim(\mathfrak L^-(V))<\infty$  if and only if   $\dim V=1$ and ${\rm ord}(p_{11})=1$ or $\infty$.
\end {Proposition}

\noindent{\it Proof.}
  Sufficiency is clear.

 Necessity. Assume $\dim V>1.$ By     Theorem \ref {5.3},  $ \Delta (\mathfrak B(V))$ is an arithmetic root system.  Consequently,  there exists $u \in D$ such that ${\rm ord }(p_{uu} )=1 $ or $\infty$. Considering Lemma \ref {3.1} (v) we get a contradiction. \hfill $\Box$

\section {Finiteness of  Nichols algebras and Nichols (braided) Lie algebras}

In this section we obtain an explicit necessary and sufficient condition for a Nichols algebra and a  Nichols (braided)  Lie algebra   to be  finite dimensional
\begin {Proposition} \label {2.4'} Assume that $\mathfrak B(V) $ is a connected  Nichols algebra of diagonal type. Then $\dim(\mathfrak B(V))=\infty$ and $\dim(\mathfrak L(V))<\infty$  if and only if   $\dim V=1$ and ${\rm ord}(p_{11})=1$.

\end {Proposition}

\noindent {\it Proof.}   Sufficiency. It is clear. Necessity.   By \cite [Th. 2.2]{WZZ15b}, $\dim V =1$. Considering $\dim(\mathfrak B(V))=\infty$ we have ${\rm ord}(p_{11})=1$ or $\infty.$  By \cite [Lemma 4.3]{WZZ15a} ${\rm ord}(p_{11})=1$.
 \hfill $\Box$

\begin {Theorem} \label {4.1} If $\mathfrak B(V) $ is connected Nichols algebra of diagonal type with $\dim V>1$,  Then $\dim(\mathfrak B(V))=\infty$  if and only if  $\dim(\mathfrak L^-(V))=\infty.$
\end {Theorem}

\noindent{\it Proof.}   Sufficiency is clear. Necessity. It follows from  Proposition \ref {3.3}. \hfill $\Box$

\vskip.1in
The fixed parameters     in \cite [Table A.1, A.2]{He05},  \cite [Table B, C]{He06a} are called  quantum numbers of  generalied Dynkin diagrams. For example,     quantum number  of Row 2 in \cite [Table A.1]{He05} is $q$; quantum numbers   of Row 9 in \cite [Table A.2]{He05} are  $q$,  $r$ and $s$.

\begin {Theorem} \label {4.2} Assume that $\mathfrak B(V) $ is connected Nichols algebra of diagonal type with $\dim V>1$,  Then the following conditions are equivalent:
{\rm (i)} $\mathfrak B(V)$ is finite-dimensional;
{\rm (ii)} $\mathfrak L(V)$ is finite-dimensional;
{\rm (iii)} $\mathfrak L^-(V)$ is finite-dimensional.
{\rm (iv)} $\Delta(\mathfrak B(V)) $ is an arithmetic root system,  ${\rm ord}(p_{ii})\not=\infty$ and ${\rm ord}(p_{ij}p_{ji})\not=\infty$ for $\forall\ 1\leq i\not=j\leq n$. {\rm (v)} $\Delta(\mathfrak B(V)) $ is an arithmetic root system and the quantum numbers (i.e. The fixed parameters  ) of generalized Dynkin diagrams of $V$ are of finite order.

\end {Theorem}
\noindent{\it Proof.}  $ {\rm (i)} \Leftrightarrow  {\rm (ii)}$. It follows from  \cite [Th. 2.2] {WZZ15b}.

$ {\rm (i)} \Leftrightarrow  {\rm (iii)}$. It follows from Theorem \ref {4.1}.

$ {\rm (iv)} \Rightarrow  {\rm (i)}$.  By Proposition  \ref {2.4'},  $p_{u,  u} \not=1$ for any $u\in D.$ We complete the proof  by  \cite [Th.1.4.6]{He05} or \cite [Th. 2] {Kh99}    since ${\rm ord }(p_{u,  u}) < \infty$ for any $u\in D$.

$ {\rm (i)} \Rightarrow  {\rm (iv)}$. Obviously,  $\Delta(\mathfrak B(V)) $ is an arithmetic root system and ${\rm ord}(p_{ii})\not=\infty$. Assume  ${\rm ord}(p_{ij}p_{ji})=\infty$ for any $i.$
By \cite [Lemma 2.2.1] {He05},  one has  $[x_i,  x_j]\in D$ or $[x_j,  x_i]\in D$. However,   the heights of $[x_i,  x_j]$ and $[x_j,  x_i]$ are infinite. Consequently we get a contradiction.

${\rm (iv)} \Leftrightarrow  {\rm (v)}$. It follows from  \cite [Table A.1,  A.2]{He05} and  \cite [Table B,  C]{He06a}. \hfill $\Box$

\begin {Corollary} \label {4.3} Assume that $\mathfrak B(V) $ is connected Nichols algebra of diagonal type with $\dim V>1$. Then  $ \dim (\mathfrak B(V)) = \infty $
(resp. $ \dim $($\mathfrak L(V)) = \infty $$)$ $($ resp. $ \dim (\mathfrak L^-(V)) = \infty $ $)$
 if and only if one of following conditions holds:

{\rm (i)} $\Delta (\mathfrak B(V)) $ is not an  arithmetic root system.

{\rm (ii)} $\dim V =2$,  ${\rm ord } (q) = \infty  $ in Row 2-6,  Row 11 of \cite [Table A.1]{He05}.

{\rm (iii)}  $\dim V =3$,  ${\rm ord } (q) = \infty  $ in Row 1-8,    10 of \cite [Table A.2]{He05}.

{\rm (iv)}  $\dim V =3$,  ${\rm ord } (q) = \infty  $ or  ${\rm ord } (r) = \infty  $ or ${\rm ord } (s) = \infty  $  in Row 9 of \cite [Table A.2]{He05}.

{\rm (v)}  $\dim V =4$,  ${\rm ord } (q) = \infty  $ in Row 1-14  of \cite [Table B]{He06a}.

{\rm (vi)}  $\dim V \ge 5$,  ${\rm ord } (q) = \infty  $ in Row 1-4,  7-10,  16,  20,  22 of \cite [Table C]{He06a}.

\end {Corollary}

\noindent{\it Proof.}  It follows from Theorem \ref {4.2},     \cite [Table A.1,  A.2]{He05} and  \cite [Table B,  C]{He06a}.\hfill $\Box$

\section{ $\mathfrak L^ - (V)$  with $\dim V=2$ and $p_{11}=p_{22}=-1$}

In this section we obtain a basis of $\mathfrak L^{-}(V)$  where  $V$ is a connected braided vector of  diagonal type with $\dim V=2$ and $p_{11}=p_{22}= -1$.

\begin {Lemma}\label {6.3}   Assume that $\mathfrak L^-(V) $ is a connected Nichols Lie algebra of diagonal type with $\dim V=2$ and $p_{11}=p_{22}= -1$.

{\rm (i)} If $u$ is a non-zero monomial, then there exists a non-zero $\alpha \in F$ such that $u =  \alpha (x_1x_2)^k x_1, $ or $\alpha x _2(x_1x_2)^k ,$
 or $\alpha (x_1x_2)^{k+1},$ or $\alpha (x_2x_1)^{k+1},$ $k\ge 0.$

{\rm (ii)} $\bar r_{x_1x_2}^i[x_1]^- = (-2)^i (x_1x_2)^ix_1$ and $\bar r_{x_2x_1}^i[x_2] = (-2)^i (x_2x_1)^ix_2$ for $i \ge 0.$

{\rm (iii)} For any method $\sigma$ adding bracket  $[, ]^-$ and  $k >0$, there exists $\alpha_k, \beta _k \in F$
such that
\begin   {eqnarray} \label {e6.3.1} \sigma (x_{i_1}, x_{i_2}, \cdots , x_{i_{2k+1}}) = \beta_k  (x_1x_2)^ {k} x_1  \hbox { or }  \beta_k  x_2(x_1x_2)^ {k}.\end   {eqnarray}
 and \begin   {eqnarray} \label {e6.3.2} \sigma (x_{i_1}, x_{i_2}, \cdots , x_{i_{2k}}) = \alpha_k  ( (x_1x_2)^ k - (x_2x_1)^ k)\end  {eqnarray}

{\rm (iv)} $$<y_1, (x_{1}x_2)^{i}-(x_{2}x_1)^{i}>=(1+(-p_{12})^{-i})x_2(x_{1}    x_2)^{i-1},$$
$$<y_2, (x_{1}x_2)^{i}-(x_{2}x_1)^{i}>=-(1+(-p_{21})^{-i})(x_{1}    x_2)^{i-1}x_1.$$

 {\rm (v)} Let  $P :=  \begin{array}{ll}
\{ x_2(x_{1}    x_2)^{k},(x_{1}  x_2)^{k}x_1, (x_{1}    x_2)^{k},  x_{2}( x_1   x_2)^{k}x_1, 0 \le k <  {\rm ord}
(p_{12} p_{21})\}.
\end{array} $ Then $P$ is a basis of $\mathfrak B(V)$.

{\rm (vi)} \begin   {eqnarray} \label {e6.3.3}
[[ [[x_{1}, x_{2}]^-, x_1 ]^-, x_2]^-\cdots , x_{1}]^-]^- = 2^{k-1} ( (x_1x_2)^k- (x_2x_1)^k   ),
\end   {eqnarray}  where there exists $2k-1$ brackets in the left hand side.
\end {Lemma}
\noindent {\it Proof.}
{\rm (i)} It is obvious.

{\rm (ii)} It can be proved by induction on $i.$

{\rm (iii)} We show this by induction on $k$. It is clear for $k =1.$ Now assume $k >1$.  Then we have $\sigma (x_{i_1}, x_{i_2}, \cdots , x_{i_{2k}})= [\sigma _1(x_{i_1}, x_{i_2}, \cdots , x_{i_{2s}}),$ $\sigma _2(x_{i_1}, x_{i_2}, \cdots , x_{i_{2t}})]^-$ or $\sigma (x_{i_1}, x_{i_2}, \cdots , x_{i_{2k}})$
$= [\sigma _1(x_{i_1}, x_{i_2}, \cdots , x_{i_{2s-1}}),$ $\sigma _2(x_{i_1}, x_{i_2}, \cdots , x_{i_{2t-1}})]^-$. In the first case, we have
\begin {eqnarray*}
\sigma (x_{i_1}, x_{i_2}, \cdots , x_{i_{2k}}) &=& [\sigma _1(x_{i_1}, x_{i_2}, \cdots , x_{i_{2s}}), \sigma _2(x_{i_{2s+1}}, x_{i_2}, \cdots , x_{i_{2k}})]^-\\
&=& [ \alpha_s ( (x_1x_2)^s-  (x_2x_1)^s), \alpha_t ( (x_1x_2)^t-  (x_2x_1)^t) ]^- \\
&& \ \ \ (\hbox {by inductive assumption, where } s +t =k)\\
&=& \alpha_k ( (x_1x_2)^k-  (x_2x_1)^k)\ \ \ (\hbox {where  } \alpha_k=0 ).
\end  {eqnarray*}

\noindent
In second case, we have
\begin {eqnarray*}
\sigma (x_{i_1}, x_{i_2}, \cdots , x_{i_{2k}}) &=& [\sigma _1(x_{i_1}, x_{i_2}, \cdots , x_{i_{2s+1}}), \sigma _2(x_{i_{2s+2}}, x_{i_2}, \cdots , x_{i_{2k}})]^-\\
&=& [ \beta _s ( (x_{j_1}x_{j_2})^s x_{j_1}, \beta _t (x_{k_1}x_{k_2})^t x_{k_1})]^- \\
&& (\hbox {by inductive assumption, where } s + t+1 =k,  j_1 \not=j_2 \hbox { and } k_1\not= k_2 ) \\
&=& \beta_s \beta _t( (x_1x_2)^k-  (x_2x_1)^k).
\end  {eqnarray*}
\begin {eqnarray*}
\sigma (x_{i_1}, x_{i_2}, \cdots , x_{i_{2k+1}}) &=& [\sigma _1(x_{i_1}, x_{i_2}, \cdots , x_{i_{2s}}), \sigma _2(x_{i_{2s+2}}, x_{i_{2s+3}}, \cdots , x_{i_{2k+1}})]^-\\
&=& [ \alpha _s ( (x_{1}x_{2})^s  - (x_{2}x_{1})^s ), \beta _t (x_{k_1}x_{k_2})^t x_{k_1})]^- \\
&& \ \ \ (\hbox {by inductive assumption, where }  s+t =k \hbox { and } k_1\not= k_2 ) \\
&=& \beta_k (x_{k_1}x_{k_2})^kx_{k_1}.
\end  {eqnarray*}
 Consequently, {\rm (iii)} holds.

{\rm (iv)}

$\begin{array}{ll}
<y_1, (x_{1}x_2)^{i}-(x_{2}x_1)^{i}>&=<y_1, x_{1}x_2(x_{1}x_2)^{i-1}
-(x_{2}x_1)^{i-1}x_{2}x_1>\\
  &=(1+(-p_{12})^{-i})x_2(x_{1}    x_2)^{i-1}.
\end {array}$

$\begin{array}{ll}
<y_2, (x_{1}x_2)^{i}-(x_{2}x_1)^{i}>&=<y_2, (x_{1}x_2)^{i-1}x_{1}x_2
-x_{2}x_1(x_{2}x_1)^{i-1}>\\
&=-(1+(-p_{21})^{-i})(x_{1}    x_2)^{i-1}x_1.
\end {array}$

{\rm (v)} It follows from \cite[Corollary 1]{Kh99}.

{\rm (vi)} We show this  by induction on $k$.
\begin {eqnarray*}
[[ [[x_{1}, x_{2}]^-, x_1 ]^-, x_2]^-\cdots , x_{1}]^-]^- &=&  [[[ 2^{k-2} ( (x_1x_2)^{k-1}-  (x_2x_1)^{k-1}), x_1]^-, x_2 ]^- \\
&& \ \ \ (\hbox {by inductive assumption})\\
&=& 2^{k-1} ( (x_1x_2)^k-  (x_2x_1)^k).
\end  {eqnarray*} \hfill $\Box$

\begin {Theorem}  \label {5.1}   Assume that $\mathfrak L^-(V) $ is connected Nichols Lie algebra of diagonal type with $\dim V=2$ and $p_{11}=p_{22}= -1$. Let\\ $P_{2}:=\left \{  \begin{array}{ll}
\{(x_{1}x_2)^{k+1}-(x_{2}x_1)^{k+1},x_2(x_{1}    x_2)^{k},(x_{1}  x_2)^{k}x_1,  0 \le k <  m \}-\{(x_{1}x_2)^{m}-(x_{2}x_1)^{m }\}, \\   \hfill  \hbox {when }
   {(-p_{12}})^{m}= {(-p_{21}})^{m}=-1 \hbox { with } {\rm ord}
(p_{12} p_{21})  =m < \infty \\
\{(x_{1}x_2)^{k+1}-(x_{2}x_1)^{k+1},x_2(x_{1}    x_2)^{k},(x_{1}  x_2)^{k}x_1,  0 \le k <  {\rm ord}
(p_{12} p_{21})  \}, \hfill  \hbox {otherwise. } \\
\end{array}\right. $
Then $P_{2}$ is a basis of $\mathfrak L^{-}(V)$. Furthermore, if  $ {\rm ord}
(p_{12} p_{21}) = m < \infty$, then

\noindent $\dim \mathfrak L^{-}(V) =\left \{  \begin{array}{ll}
3m-1,  &    \hfill  \hbox {when }
   {(-p_{12}})^{m}= {(-p_{21}})^{m}=-1,  \\
3m,   & \mbox {otherwise. } \\
\end{array}\right. .$

\end {Theorem}
\noindent {\it Proof.} By Lemma \ref {6.3} (iv)(v), $P_2$ is linearly independent.
It follows from  Lemma \ref {6.3} (ii) (vi) that $P_2 \subseteq   \mathfrak L^{-}(V).$
By Lemma \ref {6.3} (iii),  $P_{2}$ is a basis of $\mathfrak L^{-}(V)$.\hfill $\Box$

For example, if  ${\rm ord}
(p_{12} p_{21}) = m $   is a even with  $p_{12}, p_{21}\in R_{2m}$ or $m $  is a odd with $ p_{12},p_{21}\in R_{m}$, then $(-p_{12} ^m) =-1$ and $(-p_{21} ^m) =-1$. In this case, $(x_{1}x_2)^{m}-(x_{2}x_1)^{m}=0$.

\section{ Conditions for  $\overline{l}_{x_{i}}^{m}[x_{j}]^ -=0$.}
In this section we give the sufficient and necessary conditions for $\overline{l}_{x_{i}}^{m}[x_{j}]^ -=0$.

We only consider $x_1$ and $x_2$ for convenience in this section.

It is clear $\overline{l}_{x_{2}}^{m}[x_{1}]^ -=(-1)^{m}\overline{r}_{x_{2}}^{m}[x_{1}]^ -$.  We know $m_{12}\leq {\rm ord }(p_{11})-1$,
$m_{21}\leq {\rm ord }(p_{22})-1$.
 Let $i$ denote $x_{i}$ in short, sometimes.

\begin {Lemma} \label {1.0}
{\rm (i)} $x_{1}^{s}x_{2}x_{1}^{t}$ is linearly independent for $ 0\leq s\leq m_{12},0\leq t\leq{\rm ord }(p_{11})-1$.
{\rm (ii)}  $x_{2}^{s}x_{1}x_{2}^{t}$ is linearly independent for $ 0\leq t\leq m_{21},0\leq s\leq{\rm ord }(p_{22})-1$.
\end {Lemma}
\noindent {\it Proof.} {\rm (i)} It is clear
since $x_{1}^{s}x_{2}x_{1}^{t}$ is a basic element by  \cite [Lemm 2.2.1]{He05} and  \cite [Coroll 1]{Kh99}. {\rm (ii)} is similar.
\hfill $\Box$

\begin {Proposition} \label {1.1} {\rm (i)} Assume that $p_{12}=p_{21}^{-1}\not=1$, then
{\rm (1)} $\overline{l}_{1}^{m}[2]^-\not=0$ if and only if $ {\rm ord }(p_{11}) > m\geq0$.
{\rm (2)}  $\overline{l}_{2}^{m}[1]^-\not=0$ if and only if  $ {\rm ord }(p_{22}) > m\geq0$.

{\rm (ii)} {\rm (1)} Assume that $m_{12}={\rm ord }(p_{11})-1,p_{12}p_{21}\neq1$, then
$\overline{l}_{1}^{m}[2]^-\not=0$ if and only if $ 0\leq m\leq m_{12}-1+{\rm ord }(p_{11})$.
{\rm (2)} Assume that $m_{21}={\rm ord }(p_{22})-1,p_{12}p_{21}\neq1$, then $\overline{l}_{2}^{m}[1]^-\not=0$ if and only if $ 0\leq m\leq m_{21}-1+{\rm ord }(p_{22})$.
\end {Proposition}
\noindent {\it Proof.}
{\rm (i)}  We can show $\overline{l}_{1}^{m}[2]^-=(p_{12}-1)^{m}21^{m}$
and $\overline{r}_{2}^{m}[1]^-=(p_{12}-1)^{m}2^{m}1$ by induction on $m$ since $12=p_{12}21$. the other is clear by \cite [Lemma 1.3.3]{He05}.

{\rm (ii)} We obtain $\overline{l}_{1}^{m}[2]^{-}=\sum \limits _{k=0}^{m} ( - 1)^k {\scriptsize \left(\begin{array}{cc} m\\
k \end{array}\right)} x_{1} ^{m - k}x_{2} x_{1} ^k$ by \cite [Lemm 4.1(iii)]{WZZ15b} and
$\overline{r}_{2}^{m}[1]^{-}=\sum \limits _{k=0}^{m} ( - 1)^k {\scriptsize\left(\begin{array}{cc} m\\
k \end{array}\right)} x_{2} ^{k}x_{1} x_{2} ^{m-k}$ by induction on $m$.
Then $\overline{l}_{1}^{m_{12}-1+{\rm ord }(p_{11})}[2]^{-}$

\noindent $=\overline{l}_{1}^{2{\rm ord }(p_{11})-2}[2]^{-}=\sum \limits _{k=0}^{2{\rm ord }(p_{11})-2} ( - 1)^k $
${\scriptsize\left(\begin{array}{cc} 2{\rm ord }(p_{11})-2\\
k \end{array}\right)} x_i ^{2{\rm ord }(p_{11})-2 - k}x_j x_i ^k= ( - 1)^{{\rm ord }(p_{11})-1} $

\noindent $ {\scriptsize\left(\begin{array}{cc} 2{\rm ord }(p_{11})-2\\
{\rm ord }(p_{11})-1 \end{array}\right)}x_{1} ^{{\rm ord }(p_{11})-1}x_{2} x_{1} ^{{\rm ord }(p_{11})-1}$
and
$\overline{r}_{2}^{m_{21}-1+{\rm ord }(p_{22})}[1]^{-}= ( - 1)^{{\rm ord }(p_{22})-1} {\scriptsize\left(\begin{array}{cc} 2{\rm ord }(p_{22})-2\\
{\rm ord }(p_{22})-1 \end{array}\right)}$

\noindent $ x_{2} ^{{\rm ord }(p_{22})-1}x_{1} x_{2} ^{{\rm ord }(p_{22})-1}$
by \cite [Lemm 1.3.3(i)]{He05}. It is clear $\overline{l}_{1}^{m_{12}-1+{\rm ord }(p_{11})}[2]^{-}$
$\neq 0$ since $x_{1} ^{{\rm ord }(p_{11})-1}x_{2} x_{1} ^{{\rm ord }(p_{11})-1}=x_{1} ^{m_{12}}x_{2} x_{1} ^{{\rm ord }(p_{11})-1}$ is a basic element
and $\overline{r}_{2}^{m_{21}-1+{\rm ord }(p_{22})}[1]^{-}\neq 0$ since $x_{2} ^{{\rm ord }(p_{22})-1}x_{1} x_{2} ^{{\rm ord }(p_{22})-1}=x_{2} ^{{\rm ord }(p_{22})-1}x_{1} x_{2} ^{m_{21}}$ is a basic element by  Lemma \ref {1.0}. Then $\overline{l}_{1}^{m}[2]^{-}\neq 0$ if $0\leq m\leq m_{12}-1+{\rm ord }(p_{11})$ and $\overline{r}_{2}^{m}[1]^{-}\neq 0$ if $0\leq m\leq m_{21}-1+{\rm ord }(p_{22})$. On
the other hand, $\overline{l}_{1}^{m_{12}+{\rm ord }(p_{11})}[2]^{-}$
$=
( - 1)^{{\rm ord }(p_{11})-1} {\scriptsize\left(\begin{array}{cc} 2{\rm ord }(p_{11})-2\\
{\rm ord }(p_{11})-1 \end{array}\right)}[x_{1},x_{1} ^{{\rm ord }(p_{11})-1}x_{2} x_{1} ^{{\rm ord }(p_{11})-1}]=0$, $\overline{r}_{2}^{m_{21}+{\rm ord }(p_{22})}[1]^{-}$
$=
( - 1)^{{\rm ord }(p_{22})-1}$
${\scriptsize\left(\begin{array}{cc} 2{\rm ord }(p_{22})-2\\
{\rm ord }(p_{22})-1 \end{array}\right) }[x_{2} ^{{\rm ord }(p_{22})-1}x_{1} x_{2} ^{{\rm ord }(p_{22})-1},x_{2}]=0$.
\hfill $\Box$

Now we always assumed that $p_{12}p_{21}\neq1$ and $m_{12}<{\rm ord }(p_{11})-1$, i.e. $p_{11}^{m_{12}}p_{12}p_{21}=1$, $m_{21}<{\rm ord }(p_{22})-1$, i.e. $p_{22}^{m_{21}}p_{12}p_{21}=1$.

For our study of Nichols algebras we will need some non-standard formulas for quantum integers and Gaussian binomial coefficients.

In the ring $\mathbb Z[a]$, let $(0)_{a}=0$ and for any $m\in \mathbb N$, let
$(m)_{a}=1+a+a^{2} +\cdots+a^{m-1}$. The polynomials $(m)_{a}$ with $m\in \mathbb Z$ are also
known as quantum integers. Moreover, let $(0)_{a}^{!}=1$, and for any $m\in \mathbb Z$
let $(m)_{a}^{!}=\prod \limits _{i=1}^{m}(i)_{a}$. For any $i,m\in \mathbb Z$ with $0\leq i\leq m$, the rational function
${\scriptsize\left(\begin{array}{cc} m\\
k \end{array}\right)_{a}}=\frac{(m)_{a}^{!}}{(k)_{a}^{!}(m-k)_{a}^{!}},$
is in fact an element of $\mathbb Z[a]$ and is called a Gaussian binomial coefficient.
For $m\in\mathbb N_0$,$i\in \mathbb Z$ with $i<0$ or $i>m$ one defines
${\scriptsize\left(\begin{array}{cc} m\\
i \end{array}\right)_{a}}=0$. The
Gaussian binomial coefficients satisfy the following formulas:
\begin {eqnarray}\label {e1}{\scriptsize\left(\begin{array}{cc} m\\
k \end{array}\right)_{a}}={\scriptsize\left(\begin{array}{cc} m\\
m-k \end{array}\right)_{a}},\end {eqnarray}
\begin {eqnarray}\label {e2}a^{k-1}{\scriptsize\left(\begin{array}{cc} m\\
k \end{array}\right)_{a}}+a^{m}{\scriptsize\left(\begin{array}{cc} m\\
k-1 \end{array}\right)_{a}}=a^{k-1}{\scriptsize\left(\begin{array}{cc} m+1\\
k \end{array}\right)_{a}}\end {eqnarray}
\noindent for $m\in\mathbb N$, $1\leq k\leq m$.

Let $m\in\mathbb N$.
Let $f_{i,x}^{j}:=(-1)^{j}p_{11}^{\frac{j(j-1)}{2}} x^{j}{\scriptsize\left(\begin{array}{cc} i\\
j \end{array}\right)_{p_{11}}},f_{i,x}^{0}=1$.  $a_{t}:=f_{({\rm ord }(p_{11})-1)-l+t,p_{12}}^{(k-m+{\rm ord }(p_{11})-1)-l+t}-\sum \limits _{j=1}^{t-1}f_{({\rm ord }(p_{11})-1)-l+t,p_{12}}^{t-j}a_{j}$
for $t\geq0$.
 Set $g_{0,k}^{m}:=(-1)^{k} {\scriptsize\left(\begin{array}{cc} m\\
k \end{array}\right)},$
$g_{l,k}^{m}:=g_{l-1,k}^{m}
-g_{l-1,(m-{\rm ord }(p_{11})+1)+(l-1)}^{m}$

\noindent $f_{({\rm ord }(p_{11})-1)-(l-1),p_{12}}^{(k-m+{\rm ord }(p_{11})-1)-(l-1)}$ for $l\geq1$.
Let $\bar{f}_{i,x}^{j}:=(-1)^{j}p_{22}^{\frac{j(j-1)}{2}} x^{j}{\scriptsize\left(\begin{array}{cc} i\\
j \end{array}\right)_{p_{22}}},\bar{f}_{i,x}^{0}=1$.
$\bar{a}_{t}:=\bar{f}_{({\rm ord }(p_{22})-1)-l+t,p_{21}}^{(k-m+{\rm ord }(p_{22})-1)-l+t}$
$-\sum \limits _{j=1}^{t-1}\bar{f}_{({\rm ord }(p_{22})-1)-l+t,p_{21}}^{t-j}\bar{a}_{j}$
for $t\geq0$.
 Set $\bar{g}_{0,k}^{m}:=(-1)^{k} {\scriptsize\left(\begin{array}{cc} m\\
k \end{array}\right)},$
$\bar{g}_{l,k}^{m}:=\bar{g}_{l-1,k}^{m}
-\bar{g}_{l-1,(m-{\rm ord }(p_{22})+1)+(l-1)}^{m}$
$\bar{f}_{({\rm ord }(p_{22})-1)-(l-1),p_{21}}^{(k-m+{\rm ord }(p_{22})-1)-(l-1)}$ for $l\geq1$.

\begin{Lemma}\label{1.2}
{\rm (i)} $\sum \limits _{l=0}^{m}(-1)^{l}a^{\frac{l(l-1)}{2}}x^{l}{\scriptsize\left(\begin{array}{cc} m\\
l \end{array}\right)_{a}}=\prod \limits _{l=0}^{m-1}(1-a^{l}x)$
for $m\in\mathbb N$.

{\rm (ii)} $\sum \limits _{j=1}^{t}(-1)^{j-1}a^{\frac{(2-2t)(m)-2tj+j^{2}+j}{2}}{\scriptsize\left(\begin{array}{cc} m+t\\
t-j \end{array}\right)_{a}}{\scriptsize\left(\begin{array}{cc} m+j-1\\
j-1 \end{array}\right)_{a}}=1$ for $m\in\mathbb N_0$.
\end {Lemma}
\noindent {\it Proof.} {\rm (i)}  We show this by induction on $m$. If $m=1$, $1-x^{1}{\scriptsize\left(\begin{array}{cc} 1\\
1 \end{array}\right)_{a}}=1-x$, it is clear. Assume $m>1.$

$\sum \limits _{l=0}^{m}(-1)^{l}a^{\frac{l(l-1)}{2}}x^{l}{\scriptsize\left(\begin{array}{cc} m\\
l \end{array}\right)_{a}}$

$=1+\sum \limits _{l=1}^{m-1}(-1)^{l}a^{\frac{l(l-1)}{2}}x^{l}{\scriptsize\left(\begin{array}{cc} m\\
l \end{array}\right)_{a}}+(-1)^{m}a^{\frac{m(m-1)}{2}}x^{m}$

$=1+\sum \limits _{l=1}^{m-1}(-1)^{l}a^{\frac{l(l-1)}{2}}x^{l}({\scriptsize\left(\begin{array}{cc} m-1\\
l \end{array}\right)_{a}}+a^{m-l}{\scriptsize\left(\begin{array}{cc} m-1\\
l-1 \end{array}\right)_{a}})+(-1)^{m}a^{\frac{m(m-1)}{2}}x^{m}$

$=\prod \limits _{l=0}^{m-2}(1-a^{l}x)-a^{m-1}x\sum \limits _{l=0}^{m-2}(-1)^{l}a^{\frac{l(l-1)}{2}}x^{l}{\scriptsize\left(\begin{array}{cc} m-1\\
l \end{array}\right)_{a}}+(-1)^{m}a^{\frac{m(m-1)}{2}}x^{m}$

$=\prod \limits _{l=0}^{m-2}(1-a^{l}x)-a^{m-1}x(\prod \limits _{l=0}^{m-2}(1-a^{l}x)
-(-1)^{m-1}a^{\frac{(m-2)(m-1)}{2}}x^{m-1})+(-1)^{m}a^{\frac{m(m-1)}{2}}x^{m}$

$=\prod \limits _{l=0}^{m-1}(1-a^{l}x)$.

{\rm (ii)}  We show this by induction on $t$. If $t=1$, ${\scriptsize\left(\begin{array}{cc} m+1\\
0 \end{array}\right)_{a}}{\scriptsize\left(\begin{array}{cc} m\\
0\end{array}\right)_{a}}=1$, it is clear. Assume $t>1.$

$\sum \limits _{j=1}^{t}(-1)^{j-1}a^{\frac{(2-2t)(m)-2tj+j^{2}+j}{2}}{\scriptsize\left(\begin{array}{cc} m+t\\
t-j \end{array}\right)_{a}}{\scriptsize\left(\begin{array}{cc} m+j-1\\
j-1 \end{array}\right)_{a}}$

\hfill $-\sum \limits _{j=1}^{t-1}(-1)^{j-1}a^{\frac{(2-2t+2)(m)-2tj+2j+j^{2}+j}{2}}
{\scriptsize\left(\begin{array}{cc} m+t-1\\
t-1-j \end{array}\right)_{a}}{\scriptsize\left(\begin{array}{cc} m+j-1\\
j-1 \end{array}\right)_{a}}$

$=(-1)^{t-1}a^{\frac{(2-2t)(m)-t^{2}+t}{2}}{\scriptsize\left(\begin{array}{cc} m+t-1\\
t-1 \end{array}\right)_{a}}$

\hfill $+\sum \limits _{j=1}^{t-1}(-1)^{j-1}a^{\frac{(2-2t)(m)-2tj+j^{2}+j}{2}}{\scriptsize\left(\begin{array}{cc} m+t\\
t-j \end{array}\right)_{a}}{\scriptsize\left(\begin{array}{cc} m+j-1\\
j-1 \end{array}\right)_{a}}\{1-a^{m+j}\frac{{\scriptsize\left(\begin{array}{cc} m+t-1\\
t-1-j \end{array}\right)_{a}}}{{\scriptsize\left(\begin{array}{cc} m+t\\
t-j \end{array}\right)_{a}}}\}$

$=(-1)^{t-1}a^{\frac{(2-2t)(m)-t^{2}+t}{2}}{\scriptsize\left(\begin{array}{cc} m+t-1\\
t-1 \end{array}\right)_{a}}$

\hfill $+\sum \limits _{j=1}^{t-1}(-1)^{j-1}a^{\frac{(2-2t)(m)-2tj+j^{2}+j}{2}}{\scriptsize\left(\begin{array}{cc} m+t\\
t-j \end{array}\right)_{a}}{\scriptsize\left(\begin{array}{cc} m+j-1\\
j-1 \end{array}\right)_{a}}\{1-a^{m+j}\frac{(t-j)_{a}}{(m+t)_{a}}\}$

$=(-1)^{t-1}a^{\frac{(2-2t)(m)-t^{2}+t}{2}}{\scriptsize\left(\begin{array}{cc} m+t-1\\
t-1 \end{array}\right)_{a}}$

\hfill $+\sum \limits _{j=1}^{t-1}(-1)^{j-1}a^{\frac{(2-2t)(m)-2tj+j^{2}+j}{2}}{\scriptsize\left(\begin{array}{cc} m+t\\
t-j \end{array}\right)_{a}}{\scriptsize\left(\begin{array}{cc} m+j-1\\
j-1 \end{array}\right)_{a}}\frac{(m+j)_{a}}{(m+t)_{a}}$

$=\sum \limits _{j=1}^{t}(-1)^{j-1}a^{\frac{(2-2t)(m)-2tj+j^{2}+j}{2}}{\scriptsize\left(\begin{array}{cc} m+t-1\\
t-1 \end{array}\right)_{a}}{\scriptsize\left(\begin{array}{cc} t-1\\
j-1 \end{array}\right)_{a}}$

$=a^{(1-t)(m)}{\scriptsize\left(\begin{array}{cc} m+t-1\\
t-1 \end{array}\right)_{a}}\sum \limits _{j=1}^{t}(-1)^{j-1}a^{\frac{-2tj+j^{2}+j}{2}}{\scriptsize\left(\begin{array}{cc} t-1\\
j-1 \end{array}\right)_{a}}$

$=a^{(1-t)(m)}{\scriptsize\left(\begin{array}{cc} m+t-1\\
t-1 \end{array}\right)_{a}}\sum \limits _{j=0}^{t-1}(-1)^{j}a^{\frac{-2t(j+1)+(j+1)^{2}+(j+1)}{2}}{\scriptsize\left(\begin{array}{cc} t-1\\
j \end{array}\right)_{a}}$

$=a^{(1-t)(m)}{\scriptsize\left(\begin{array}{cc} m+t-1\\
t-1 \end{array}\right)_{a}}a^{1-t}\sum \limits _{j=0}^{t-1}(-1)^{j}a^{\frac{j(j-1)}{2}}
(a^{2-t})^{j}{\scriptsize\left(\begin{array}{cc} t-1\\
j \end{array}\right)_{a}}$

$=a^{(1-t)(m+1)}{\scriptsize\left(\begin{array}{cc} m+t-1\\
t-1 \end{array}\right)_{a}}\prod \limits _{j=0}^{t-2}(1-a^{j}a^{2-t})=0$.
\hfill $\Box$

 Set $b_{1,k}^{m}:=f_{m,p_{12}}^{k},$
$b_{l,k}^{m}:=b_{l-1,k}^{m}
-b_{l-1,l-1}^{m}f_{m-(l-1),p_{12}}^{k-(l-1)}$ for $l>1$. Let $b_{t}:=a_{t}$ for $t\geq0$ while replace ${\rm ord }(p_{11})-1$ in the definition of $a_{t}$ with $m$, i.e. $b_{t}:=f_{m-l+t,p_{12}}^{k-l+t}-\sum \limits _{j=1}^{t-1}f_{m-l+t,p_{12}}^{t-j}b_{j}$
for $t\geq0$. Set $\bar{b}_{1,k}^{m}:=\bar{f}_{m,p_{21}}^{k},$
$\bar{b}_{l,k}^{m}:=\bar{b}_{l-1,k}^{m}
-\bar{b}_{l-1,l-1}^{m}\bar{f}_{m-(l-1),p_{21}}^{k-(l-1)}$ for $l>1$. Let $\bar{b}_{t}:=\bar{a}_{t}$ for $t\geq0$ while replace ${\rm ord }(p_{22})-1$ in the definition of $\bar{a}_{t}$ with $m$, i.e. $\bar{b}_{t}:=\bar{f}_{m-l+t,p_{21}}^{k-l+t}-\sum \limits _{j=1}^{t-1}\bar{f}_{m-l+t,p_{21}}^{t-j}\bar{b}_{j}$
for $t\geq0$.

\begin {Lemma} \label {1.3}
{\rm (i)} {\rm (1)} $g_{l,k}^{m}=g_{0,k}^{m}
-\sum \limits _{j=0}^{l-1}g_{j,(m-{\rm ord }(p_{11})+1)+j}^{m}f_{({\rm ord }(p_{11})-1)-j,p_{12}}^{(k-m+{\rm ord }(p_{11})-1)-j}$. $\bar{g}_{l,k}^{m}=\bar{g}_{0,k}^{m}
-\sum \limits _{j=0}^{l-1}\bar{g}_{j,(m-{\rm ord }(p_{22})+1)+j}^{m}\bar{f}_{({\rm ord }(p_{22})-1)-j,p_{21}}^{(k-m+{\rm ord }(p_{22})-1)-j}$.

{\rm (2)} $b_{l,k}^{m}=b_{1,k}^{m}
-\sum \limits _{j=1}^{l-1}b_{j,j}^{m}f_{m-j,p_{12}}^{k-j}$.
$\bar{b}_{l,k}^{m}=\bar{b}_{1,k}^{m}
-\sum \limits _{j=1}^{l-1}\bar{b}_{j,j}^{m}\bar{f}_{m-j,p_{21}}^{k-j}$.

{\rm (ii)} {\rm (1)} $g_{l,k}^{m}=g_{0,k}^{m}
-\sum \limits _{s=1}^{r}g_{0,(m-{\rm ord }(p_{11})+1)+l-s}^{m}a_{s}-\sum \limits _{j=0}^{l-r-1}g_{j,(m-{\rm ord }(p_{11})+1)+j}^{m}
$
$(f_{({\rm ord }(p_{11})-1)-j,p_{12}}^{(k-m+{\rm ord }(p_{11})-1)-j}-\sum \limits _{s=1}^{r}f_{({\rm ord }(p_{11})-1)-j,p_{12}}^{l-s-j}a_{s})$ for $0\leq r\leq l$.
In particular,
$g_{l,k}^{m}=g_{0,k}^{m}
-\sum \limits _{j=0}^{l-1}g_{0,(m-{\rm ord }(p_{11})+1)+j}^{m}a_{l-j}$.
$\bar{g}_{l,k}^{m}=\bar{g}_{0,k}^{m}
-\sum \limits _{s=1}^{r}\bar{g}_{0,(m-{\rm ord }(p_{22})+1)+l-s}^{m}\bar{a}_{s}-\sum \limits _{j=0}^{l-r-1}\bar{g}_{j,(m-{\rm ord }(p_{22})+1)+j}^{m}
$
$(\bar{f}_{({\rm ord }(p_{22})-1)-j,p_{21}}^{(k-m+{\rm ord }(p_{22})-1)-j}$

\noindent $-\sum \limits _{s=1}^{r}\bar{f}_{({\rm ord }(p_{22})-1)-j,p_{21}}^{l-s-j}\bar{a}_{s})$ for $0\leq r\leq l$.
In particular,
$\bar{g}_{l,k}^{m}=\bar{g}_{0,k}^{m}
-\sum \limits _{j=0}^{l-1}\bar{g}_{0,(m-{\rm ord }(p_{22})+1)+j}^{m}\bar{a}_{l-j}$.

{\rm (2)} $b_{l,k}^{m}=b_{1,k}^{m}
-\sum \limits _{s=1}^{r}b_{1,l-s}^{m}b_{s}-\sum \limits _{j=1}^{l-r-1}b_{j,j}^{m}
(f_{m-j,p_{12}}^{k-j}-\sum \limits _{s=1}^{r}f_{m-j,p_{12}}^{l-s-j}b_{s})$ for $0\leq r\leq l$.
In particular, $b_{l,k}^{m}=b_{1,k}^{m}
-\sum \limits _{j=0}^{l-1}b_{1,j}^{m}b_{l-j}$.
$\bar{b}_{l,k}^{m}=\bar{b}_{1,k}^{m}
-\sum \limits _{s=1}^{r}\bar{b}_{1,l-s}^{m}\bar{b}_{s}-\sum \limits _{j=1}^{l-r-1}\bar{b}_{j,j}^{m}
(\bar{f}_{m-j,p_{21}}^{k-j}-\sum \limits _{s=1}^{r}\bar{f}_{m-j,p_{21}}^{l-s-j}\bar{b}_{s})$ for $0\leq r\leq l$.
In particular,
$\bar{b}_{l,k}^{m}=\bar{b}_{1,k}^{m}
-\sum \limits _{j=0}^{l-1}\bar{b}_{1,j}^{m}\bar{b}_{l-j}$.

{\rm (iii)} {\rm (1)} $a_{t}=f_{({\rm ord }(p_{11})-1)-l+t,p_{12}}^{(k-m+{\rm ord }(p_{11})-1)-l+t}(-1)^{t-1}p_{11}^{-(t-1)((k-m+{\rm ord }(p_{11})-1)-l)-\frac{t(t-1)}{2}}
$

${\scriptsize\left(\begin{array}{cc} (k-m+{\rm ord }(p_{11})-1)-l+t-1\\
t-1 \end{array}\right)_{p_{11}}}$, i.e.

$a_{t}=(-1)^{(k-m+{\rm ord }(p_{11})-1)-l-1}p_{11}^{\frac{((k-m+{\rm ord }(p_{11})-1)-l+1)((k-m+{\rm ord }(p_{11})-1)-l)}{2}}p_{12}^{(k-m+{\rm ord }(p_{11})-1)-l+t}$

$
{\scriptsize\left(\begin{array}{cc} ({\rm ord }(p_{11})-1)-l+t\\
(k-m+{\rm ord }(p_{11})-1)-l+t \end{array}\right)_{p_{11}}}{\scriptsize\left(\begin{array}{cc} (k-m+{\rm ord }(p_{11})-1)-l+t-1\\
t-1 \end{array}\right)_{p_{11}}}$ for $t\geq1$.

$\bar{a}_{t}=\bar{f}_{({\rm ord }(p_{22})-1)-l+t,p_{21}}^{(k-m+{\rm ord }(p_{22})-1)-l+t}(-1)^{t-1}p_{22}^{-(t-1)((k-m+{\rm ord }(p_{22})-1)-l)-\frac{t(t-1)}{2}}
$

${\scriptsize\left(\begin{array}{cc} (k-m+{\rm ord }(p_{22})-1)-l+t-1\\
t-1 \end{array}\right)_{p_{22}}}$, i.e.

$\bar{a}_{t}=(-1)^{(k-m+{\rm ord }(p_{22})-1)-l-1}p_{22}^{\frac{((k-m+{\rm ord }(p_{22})-1)-l+1)((k-m+{\rm ord }(p_{22})-1)-l)}{2}}p_{21}^{(k-m+{\rm ord }(p_{22})-1)-l+t}$

$
{\scriptsize\left(\begin{array}{cc} ({\rm ord }(p_{22})-1)-l+t\\
(k-m+{\rm ord }(p_{22})-1)-l+t \end{array}\right)_{p_{22}}}{\scriptsize\left(\begin{array}{cc} (k-m+{\rm ord }(p_{22})-1)-l+t-1\\
t-1 \end{array}\right)_{p_{22}}}$ for $t\geq1$.

 {\rm (2)} $b_{t}=(-1)^{k-l-1}p_{11}^{\frac{(k-l+1)(k-l)}{2}}p_{12}^{k-l+t}$
$
{\scriptsize\left(\begin{array}{cc} m-l+t\\
k-l+t \end{array}\right)_{p_{11}}}{\scriptsize\left(\begin{array}{cc} k-l+t-1\\
t-1 \end{array}\right)_{p_{11}}}$ for $t\geq1$.

$\bar{b}_{t}=(-1)^{k-l-1}p_{22}^{\frac{(k-l+1)(k-l)}{2}}p_{21}^{k-l+t}$
$
{\scriptsize\left(\begin{array}{cc} m-l+t\\
k-l+t \end{array}\right)_{p_{22}}}{\scriptsize\left(\begin{array}{cc} k-l+t-1\\
t-1 \end{array}\right)_{p_{22}}}$ for $t\geq1$.

{\rm (iv)} {\rm (1)} $g_{l,k}^{m}=(-1)^{k} {\scriptsize\left(\begin{array}{cc} m\\
k \end{array}\right)}
-\sum \limits _{j=0}^{l-1}(-1)^{(m-{\rm ord }(p_{11})+1)+j} {\scriptsize\left(\begin{array}{cc} m\\
(m-{\rm ord }(p_{11})+1)+j \end{array}\right)}$

$(-1)^{(k-m+{\rm ord }(p_{11})-1)-l-1}p_{11}^{\frac{((k-m+{\rm ord }(p_{11})-1)-l+1)((k-m+{\rm ord }(p_{11})-1)-l)}{2}}p_{12}^{(k-m+{\rm ord }(p_{11})-1)-j}$

$
{\scriptsize\left(\begin{array}{cc} ({\rm ord }(p_{11})-1)-j\\
(k-m+{\rm ord }(p_{11})-1)-j \end{array}\right)_{p_{11}}}{\scriptsize\left(\begin{array}{cc} (k-m+{\rm ord }(p_{11})-1)-j-1\\
l-j-1 \end{array}\right)_{p_{11}}}$.

$\bar{g}_{l,k}^{m}=(-1)^{k} {\scriptsize\left(\begin{array}{cc} m\\
k \end{array}\right)}
-\sum \limits _{j=0}^{l-1}(-1)^{(m-{\rm ord }(p_{22})+1)+j} {\scriptsize\left(\begin{array}{cc} m\\
(m-{\rm ord }(p_{22})+1)+j \end{array}\right)}$

$(-1)^{(k-m+{\rm ord }(p_{22})-1)-l-1}p_{22}^{\frac{((k-m+{\rm ord }(p_{22})-1)-l+1)((k-m+{\rm ord }(p_{22})-1)-l)}{2}}p_{21}^{(k-m+{\rm ord }(p_{22})-1)-j}$

$
{\scriptsize\left(\begin{array}{cc} ({\rm ord }(p_{22})-1)-j\\
(k-m+{\rm ord }(p_{22})-1)-j \end{array}\right)_{p_{22}}}{\scriptsize\left(\begin{array}{cc} (k-m+{\rm ord }(p_{22})-1)-j-1\\
l-j-1 \end{array}\right)_{p_{22}}}$.

{\rm (2)} $b_{l,k}^{m}=(-1)^{k}p_{11}^{\frac{k(k-1)}{2}} p_{12}^{k}{\scriptsize\left(\begin{array}{cc} m\\
k \end{array}\right)_{p_{11}}}
-\sum \limits _{j=0}^{l-1}(-1)^{j}p_{11}^{\frac{j(j-1)}{2}} p_{12}^{j}{\scriptsize\left(\begin{array}{cc} m\\
j \end{array}\right)_{p_{11}}}(-1)^{k-l-1}p_{11}^{\frac{(k-l+1)(k-l)}{2}}p_{12}^{k-j}$
$
{\scriptsize\left(\begin{array}{cc} m-j\\
k-j \end{array}\right)_{p_{11}}}$
${\scriptsize\left(\begin{array}{cc} k-j-1\\
l-j-1 \end{array}\right)_{p_{11}}}$.
$\bar{b}_{l,k}^{m}=(-1)^{k}p_{22}^{\frac{k(k-1)}{2}} p_{21}^{k}{\scriptsize\left(\begin{array}{cc} m\\
k \end{array}\right)_{p_{22}}}
-\sum \limits _{j=0}^{l-1}(-1)^{j}p_{22}^{\frac{j(j-1)}{2}} p_{21}^{j}{\scriptsize\left(\begin{array}{cc} m\\
j \end{array}\right)_{p_{22}}}$

\noindent $(-1)^{k-l-1}p_{22}^{\frac{(k-l+1)(k-l)}{2}}p_{21}^{k-j}$
$
{\scriptsize\left(\begin{array}{cc} m-j\\
k-j \end{array}\right)_{p_{22}}}$
${\scriptsize\left(\begin{array}{cc} k-j-1\\
l-j-1 \end{array}\right)_{p_{22}}}$.
\end {Lemma}
\noindent {\it Proof.}
{\rm (ii)} {\rm (1)} We show this by induction on $r$. If $r=0$, it is clear. Assume $r>0.$

$g_{l,k}^{m}=g_{0,k}^{m}
-\sum \limits _{s=1}^{r}g_{0,(m-{\rm ord }(p_{11})+1)+l-s}^{m}a_{s}-\sum \limits _{j=0}^{l-r-1}g_{j,(m-{\rm ord }(p_{11})+1)+j}^{m}
$

 \hfill $(f_{({\rm ord }(p_{11})-1)-j,p_{12}}^{(k-m+{\rm ord }(p_{11})-1)-j}-\sum \limits _{s=1}^{r}f_{({\rm ord }(p_{11})-1)-j,p_{12}}^{l-s-j}a_{s})$

$=g_{0,k}^{m}
-\sum \limits _{s=1}^{r}g_{0,(m-{\rm ord }(p_{11})+1)+l-s}^{m}a_{s}$

 \hfill $-\sum \limits _{j=0}^{l-r-2}g_{j,(m-{\rm ord }(p_{11})+1)+j}^{m}(f_{({\rm ord }(p_{11})-1)-j,p_{12}}^{(k-m+{\rm ord }(p_{11})-1)-j}
$
$-\sum \limits _{s=1}^{r}f_{({\rm ord }(p_{11})-1)-j,p_{12}}^{l-s-j}a_{s})$

 \hfill $-g_{l-r-1,(m-{\rm ord }(p_{11})+1)+l-r-1}^{m}
(f_{({\rm ord }(p_{11})-1)-l+r+1,p_{12}}^{(k-m+{\rm ord }(p_{11})-1)-l+r+1}-\sum \limits _{s=1}^{r}f_{({\rm ord }(p_{11})-1)-l+r+1,p_{12}}^{-s+r+1}a_{s})$

$=g_{0,k}^{m}
-\sum \limits _{s=1}^{r}g_{0,(m-{\rm ord }(p_{11})+1)+l-s}^{m}a_{s}-\sum \limits _{j=0}^{l-r-2}g_{j,(m-{\rm ord }(p_{11})+1)+j}^{m}
(f_{({\rm ord }(p_{11})-1)-j,p_{12}}^{(k-m+{\rm ord }(p_{11})-1)-j}$

$-\sum \limits _{s=1}^{r}f_{({\rm ord }(p_{11})-1)-j,p_{12}}^{l-s-j}a_{s})$
$-(g_{0,(m-{\rm ord }(p_{11})+1)+l-r-1}^{m}
-\sum \limits _{j=0}^{l-r-2}g_{j,(m-{\rm ord }(p_{11})+1)+j}^{m}f_{({\rm ord }(p_{11})-1)-j,p_{12}}^{l-(r+1)-j})a_{r+1}$

$=g_{0,k}^{m}
-\sum \limits _{s=1}^{r}g_{0,(m-{\rm ord }(p_{11})+1)+l-s}^{m}a_{s}-g_{0,(m-{\rm ord }(p_{11})+1)+l-r-1}^{m}a_{r+1}
-\sum \limits _{j=0}^{l-r-2}g_{j,(m-{\rm ord }(p_{11})+1)+j}^{m}
$

 \hfill $(f_{({\rm ord }(p_{11})-1)-j,p_{12}}^{(k-m+{\rm ord }(p_{11})-1)-j}-\sum \limits _{s=1}^{r}f_{({\rm ord }(p_{11})-1)-j,p_{12}}^{l-s-j}a_{s}-f_{({\rm ord }(p_{11})-1)-j,p_{12}}^{(l-r-1)-j}a_{r+1})$

$=g_{0,k}^{m}
-\sum \limits _{s=1}^{r+1}g_{0,(m-{\rm ord }(p_{11})+1)+l-s}^{m}a_{s}
-\sum \limits _{j=0}^{l-r-2}g_{j,(m-{\rm ord }(p_{11})+1)+j}^{m}
(f_{({\rm ord }(p_{11})-1)-j,p_{12}}^{(k-m+{\rm ord }(p_{11})-1)-j}$

\noindent $-\sum \limits _{s=1}^{r+1}f_{({\rm ord }(p_{11})-1)-j,p_{12}}^{l-s-j}a_{s})$, it is proved.

{\rm (2)} We show this by induction on $r$. If $r=0$, it is clear. Assume $r>0.$

$b_{l,k}^{m}=b_{1,k}^{m}
-\sum \limits _{s=1}^{r}b_{1,l-s}^{m}b_{s}-\sum \limits _{j=1}^{l-r-1}b_{j,j}^{m}
(f_{m-j,p_{12}}^{k-j}-\sum \limits _{s=1}^{r}f_{m-j,p_{12}}^{l-s-j}b_{s})$

$=b_{1,k}^{m}
-\sum \limits _{s=1}^{r}b_{1,l-s}^{m}b_{s}-\sum \limits _{j=1}^{l-r-2}b_{j,j}^{m}
(f_{m-j,p_{12}}^{k-j}-\sum \limits _{s=1}^{r}f_{m-j,p_{12}}^{l-s-j}b_{s})$

\hfill $-b_{l-r-1,l-r-1}^{m}
(f_{m-(l-r-1),p_{12}}^{k-(l-r-1)}-\sum \limits _{s=1}^{r}f_{m-(l-r-1),p_{12}}^{l-s-(l-r-1)}b_{s})$

$=b_{1,k}^{m}
-\sum \limits _{s=1}^{r}b_{1,l-s}^{m}b_{s}-\sum \limits _{j=1}^{l-r-2}b_{j,j}^{m}
(f_{m-j,p_{12}}^{k-j}-\sum \limits _{s=1}^{r}f_{m-j,p_{12}}^{l-s-j}b_{s})$

\hfill $-(b_{1,l-r-1}^{m}
-\sum \limits _{j=1}^{l-r-2}b_{j,j}^{m}f_{m-j,p_{12}}^{l-r-1-j})
(f_{m-l+r+1,p_{12}}^{k-l+r+1}-\sum \limits _{s=1}^{r}f_{m-l+r+1,p_{12}}^{r+1-s}b_{s})$

$=b_{1,k}^{m}
-\sum \limits _{s=1}^{r}b_{1,l-s}^{m}b_{s}-\sum \limits _{j=1}^{l-r-2}b_{j,j}^{m}
(f_{m-j,p_{12}}^{k-j}-\sum \limits _{s=1}^{r}f_{m-j,p_{12}}^{l-s-j}b_{s})-(b_{1,l-r-1}^{m}
-\sum \limits _{j=1}^{l-r-2}b_{j,j}^{m}f_{m-j,p_{12}}^{l-r-1-j})
b_{r+1}$

$=b_{1,k}^{m}
-\sum \limits _{s=1}^{r+1}b_{1,l-s}^{m}b_{s}-\sum \limits _{j=1}^{l-r-2}b_{j,j}^{m}
(f_{m-j,p_{12}}^{k-j}-\sum \limits _{s=1}^{r+1}f_{m-j,p_{12}}^{l-s-j}b_{s})$
, it is proved.

{\rm (iii)}  We show this by induction on $t$. If $t=1$, it is clear. Assume $t>1.$

$a_{t}=f_{({\rm ord }(p_{11})-1)-l+t,p_{12}}^{(k-m+{\rm ord }(p_{11})-1)-l+t}-\sum \limits _{j=1}^{t-1}f_{({\rm ord }(p_{11})-1)-l+t,p_{12}}^{t-j}a_{j}$

$=f_{({\rm ord }(p_{11})-1)-l+t,p_{12}}^{(k-m+{\rm ord }(p_{11})-1)-l+t}-\sum \limits _{j=1}^{t-1}f_{({\rm ord }(p_{11})-1)-l+t,p_{12}}^{t-j}f_{({\rm ord }(p_{11})-1)-l+j,p_{12}}^{(k-m+{\rm ord }(p_{11})-1)-l+j}(-1)^{j-1}$

\hfill $p_{11}^{-(j-1)((k-m+{\rm ord }(p_{11})-1)-l)-\frac{j(j-1)}{2}}
{\scriptsize\left(\begin{array}{cc} (k-m+{\rm ord }(p_{11})-1)-l+j-1\\
j-1 \end{array}\right)_{p_{11}}}$

$=(-1)^{(k-m+{\rm ord }(p_{11})-1)-l+t}p_{11}^{\frac{((k-m+{\rm ord }(p_{11})-1)-l+t)((k-m+{\rm ord }(p_{11})-1)-l+t-1)}{2}} p_{12}^{(k-m+{\rm ord }(p_{11})-1)-l+t}$

${\scriptsize\left(\begin{array}{cc} ({\rm ord }(p_{11})-1)-l+t\\
(k-m+{\rm ord }(p_{11})-1)-l+t \end{array}\right)_{p_{11}}}$
$-\sum \limits _{j=1}^{t-1}(-1)^{t-j}p_{11}^{\frac{(t-j)(t-j-1)}{2}} p_{12}^{t-j}{\scriptsize\left(\begin{array}{cc}({\rm ord }(p_{11})-1)-l+t\\
t-j \end{array}\right)_{p_{11}}}$

$(-1)^{(k-m+{\rm ord }(p_{11})-1)-l+j}p_{11}^{\frac{((k-m+{\rm ord }(p_{11})-1)-l+j)
((k-m+{\rm ord }(p_{11})-1)-l+j-1)}{2}} p_{12}^{(k-m+{\rm ord }(p_{11})-1)-l+j}(-1)^{j-1}$

${\scriptsize\left(\begin{array}{cc}({\rm ord }(p_{11})-1)-l+j\\
(k-m+{\rm ord }(p_{11})-1)-l+j \end{array}\right)_{p_{11}}}$
$p_{11}^{-(j-1)((k-m+{\rm ord }(p_{11})-1)-l)-\frac{j(j-1)}{2}}
{\scriptsize\left(\begin{array}{cc} (k-m+{\rm ord }(p_{11})-1)-l+j-1\\
j-1 \end{array}\right)_{p_{11}}}$

$=f_{({\rm ord }(p_{11})-1)-l+t,p_{12}}^{(k-m+{\rm ord }(p_{11})-1)-l+t}$
$\{1-\sum \limits _{j=1}^{t-1}(-1)^{j-1}p_{11}^{\frac{j(j+1)}{2}}
p_{11}^{-(t-1)(k-m+{\rm ord }(p_{11})-1-l)}
$

\hfill ${\scriptsize\left(\begin{array}{cc}(k-m+{\rm ord }(p_{11})-1)-l+t\\
t-j \end{array}\right)_{p_{11}}}$
$p_{11}^{-tj}
{\scriptsize\left(\begin{array}{cc} (k-m+{\rm ord }(p_{11})-1)-l+j-1\\
j-1 \end{array}\right)_{p_{11}}}\}$

$=f_{({\rm ord }(p_{11})-1)-l+t,p_{12}}^{(k-m+{\rm ord }(p_{11})-1)-l+t}$
$(-1)^{t-1}p_{11}^{\frac{t(t+1)}{2}-t^{2}}
p_{11}^{-(t-1)(k-m+{\rm ord }(p_{11})-1-l)}
{\scriptsize\left(\begin{array}{cc} (k-m+{\rm ord }(p_{11})-1)-l+t-1\\
t-1 \end{array}\right)_{p_{11}}}$

$=f_{({\rm ord }(p_{11})-1)-l+t,p_{12}}^{(k-m+{\rm ord }(p_{11})-1)-l+t}$
$(-1)^{t-1}p_{11}^{\frac{t(1-t)}{2}-(t-1)(k-m+{\rm ord }(p_{11})-1-l)}
{\scriptsize\left(\begin{array}{cc} (k-m+{\rm ord }(p_{11})-1)-l+t-1\\
t-1 \end{array}\right)_{p_{11}}}$.
\hfill $\Box$

\begin{Lemma}\label{6.5} Let $m\in\mathbb N$.

{\rm (i)} $f_{m,p_{12}}^{k}-p_{11}^{m}p_{12}f_{m,p_{12}}^{k-1}=f_{m+1,p_{12}}^{k}$,
$r_{1}^{m}[2]=\sum \limits _{k=0}^{m}f_{m,p_{12}}^{k}x_{1} ^{m-k}x_{2} x_{1} ^{k}$;
$\bar{f}_{m,p_{21}}^{k}-p_{22}^{m}p_{21}\bar{f}_{m,p_{21}}^{k-1}=\bar{f}_{m+1,p_{21}}^{k}$,
$r_{2}^{m}[1]=\sum \limits _{k=0}^{m}\bar{f}_{m,p_{21}}^{k}x_{2} ^{k}x_{1} x_{2} ^{m-k}$.

{\rm (ii)}
{\rm (1)} $x_{1} ^{m}x_{2} =-\sum \limits _{k=1}^{m}f_{m,p_{12}}^{k}x_{1} ^{m-k}x_{2}
 x_{1}^{k}$  for $\forall\ m>m_{12}$,
$x_{1} x_{2} ^{m}=-\sum \limits _{k=1}^{m}\bar{f}_{m,p_{21}}^{k}x_{2} ^{k}x_{1} x_{2} ^{m-k}$ for $\forall\ m>m_{21}$.

{\rm (2)} $x_{1} ^{m}x_{2} =-\sum \limits _{k=l}^{m}b_{l,k}^{m}x_{1} ^{m-k}x_{2} x_{1}^{k}$
for $\forall\ m-m_{12}\geq l\geq 1$. $x_{1}x_{2} ^{m} =-\sum \limits _{k=l}^{m}\bar{b}_{l,k}^{m}x_{2} ^{k}x_{1} x_{2} ^{m-k}$
for $\forall\ m-m_{21}\geq l\geq 1$. Specially, $x_{1} ^{m}x_{2} =-\sum \limits _{k=0}^{m_{12}}b_{m-m_{12},m-k}^{m}x_{1} ^{k}x_{2} x_{1} ^{m-k}$.
for $\forall\ m>m_{12}$. $x_{1}x_{2} ^{m} =-\sum \limits _{k=0}^{m_{21}}\bar{b}_{m-m_{21},m-k}^{m}x_{2} ^{m-k}x_{1} x_{2} ^{k}$.
for $\forall\ m>m_{21}$.

{\rm (iii)}
\noindent
$\overline{l}_{1}^{m}[2]^{-} =
\sum \limits _{k=m-{\rm ord }(p_{11})+l+1}^{{\rm ord }(p_{11})-1} g_{l,k}^{m} x_{1} ^{m- k}x_{2} x_{1} ^{k}$ for $ m\geq {\rm ord }(p_{11}),{\rm ord }(p_{11})-m_{12}-1\geq l\geq0$.

$\overline{l}_{2}^{m}[1]^{-} =(-1)^{m}
\sum \limits _{k=m-{\rm ord }(p_{22})+l+1}^{{\rm ord }(p_{22})-1} \bar{g}_{l,k}^{m} x_{2} ^{ k}x_{1} x_{2} ^{m-k}
$ for $m\geq {\rm ord }(p_{22}),{\rm ord }(p_{22})-m_{21}-1\geq l\geq0$.

{\rm (iv)}  $ \overline{l}_{1}^{m}[2]^{-}=\sum \limits _{k=m-m_{12}}^{{\rm ord }(p_{11})-1} g_{{\rm ord }(p_{11})-m_{12}-1,k}^{m} x_{1} ^{m- k}x_{2} x_{1} ^k$ for $m\geq {\rm ord }(p_{11}).$

$ \overline{l}_{2}^{m}[1]^{-}=(-1)^{m}\sum \limits _{k=m-m_{21}}^{{\rm ord }(p_{22})-1} \bar{g}_{{\rm ord }(p_{22})-m_{21}-1,k}^{m} x_{2} ^{k}x_{1} x_{2} ^{m-k}$ for $m\geq {\rm ord }(p_{22}).$

{\rm (v)} For $0\leq l<m_{12}$.
$ \overline{l}_{1}^{{\rm ord }(p_{11})+l}[2]^{-}=0$ if
and only if $g_{{\rm ord }(p_{11})-1-m_{12},{\rm ord }(p_{11})+l-m_{12}+k}^{{\rm ord }(p_{11})+l}=0$ for $\forall\ k=0,1,\ldots,m_{12}-l-1$.
For $0\leq l<m_{21}$.
$ \overline{l}_{2}^{{\rm ord }(p_{22})+l}[1]^{-}=0$ if
and only if $\bar{g}_{{\rm ord }(p_{22})-1-m_{21},{\rm ord }(p_{22})+l-m_{21}+k}^{{\rm ord }(p_{22})+l}=0$ for $\forall\ k=0,1,\ldots,m_{21}-l-1$.
\end {Lemma}
\noindent {\it Proof.}
{\rm (i)}  We show this by induction on $m$. If $m=1$, $r_{1}[2]=x_{1} x_{2}+(-1) p_{12}x_{2} x_{1} $, it is clear. Assume $m>1.$

$r_{1}^{m+1}[2]=x_{1} r_{1}^{m}[2]-p_{11}^{m}p_{12}r_{1}^{m}[2]x_{1}=\sum \limits _{k=0}^{m}f_{m,p_{12}}^{k}x_{1} ^{m-k+1}x_{2} x_{1}^{k}-p_{11}^{m}p_{12}\sum \limits _{k=0}^{m}f_{m,p_{12}}^{k}x_{1} ^{m-k}x_{2} x_{1} ^{k+1}$

$=x_{1} ^{m+1}x_{2}+\sum \limits _{k=1}^{m}f_{m,p_{12}}^{k}x_{1} ^{m-k+1}x_{2} x_{1} ^{k}-p_{11}^{m}p_{12}\sum \limits _{k=0}^{m-1}f_{m,p_{12}}^{k}x_{1} ^{m-k}x_{2} x_{1} ^{k+1}-p_{11}^{m}p_{12}f_{m,p_{12}}^{m}x_{2} x_{1} ^{m+1}$

$=x_{1} ^{m+1}x_{2}+\sum \limits _{k=1}^{m}f_{m,p_{12}}^{k}x_{1} ^{m-k+1}x_{2} x_{1} ^{k}-p_{11}^{m}p_{12}\sum \limits _{k=1}^{m}f_{m,p_{12}}^{k-1}x_{1} ^{m-k+1}x_{2} x_{1} ^{k}-p_{11}^{m}p_{12}(-1)^{m}p_{11}^{\frac{m(m-1)}{2}} p_{12}^{m}x_{2} x_{1} ^{m+1}$

$=x_{1} ^{m+1}x_{2}+\sum \limits _{k=1}^{m}f_{m+1,p_{12}}^{k}x_{1} ^{m-k+1}x_{2} x_{1} ^{k}+f_{m+1,p_{12}}^{m+1}x_{2} x_{1} ^{m+1}$

$=\sum \limits _{k=0}^{m+1}f_{m+1,p_{12}}^{k}x_{1} ^{m+1-k}x_{2} x_{1} ^{k}$.
It is proved.

{\rm (ii)} {\rm (1)} It is clear by  \cite [Lemm 4.14(ii)]{WZZ15a}.

{\rm (2)} $l=1$ is clear by {\rm (1)}, if $ m-m_{12}> l$,
$x_{1} ^{m}x_{2} =-\sum \limits _{k=l}^{m}b_{l,k}^{m}x_{1} ^{m-k}x_{2} x_{1} ^{k}=-\sum \limits _{k=l+1}^{m}b_{l,k}^{m}x_{1} ^{m-k}x_{2} x_{1} ^{k}-b_{l,l}^{m}x_{1} ^{m-l}x_{2} x_{1} ^{l}
=-\sum \limits _{k=l+1}^{m}b_{l,k}^{m}x_{1} ^{m-k}x_{2} x_{1} ^{k}
+b_{l,l}^{m}\sum \limits _{k=1}^{m-l}f_{m-l,p_{12}}^{k}x_{1} ^{m-l-k}x_{2} x_{1} ^{k}x_{1} ^{l}
=-\sum \limits _{k=l+1}^{m}b_{l,k}^{m}x_{1} ^{m-k}x_{2} x_{1} ^{k}
+b_{l,l}^{m}\sum \limits _{k=1+l}^{m}f_{m-l,p_{12}}^{k-l}x_{1} ^{m-k}x_{2} x_{1} ^{k}
=-\sum \limits _{k=l+1}^{m}(b_{l,k}^{m}
-b_{l,l}^{m}f_{m-l,p_{12}}^{k-l})x_{1} ^{m-k}x_{2} x_{1} ^{k}
=-\sum \limits _{k=l+1}^{m}b_{l+1,k}^{m}x_{1} ^{m-k}x_{2} x_{1} ^{k}$. The other is clear.

{\rm (iii)} We show this by induction on $l$ for ${\rm ord }(p_{11})-m_{12}-1\geq l\geq0$. See

$\overline{l}_{1}^{m}[2]^{-}=\sum \limits _{k=0}^{m} ( - 1)^k \left(\begin{array}{cc} m\\
k \end{array}\right) x_{1} ^{m- k}x_{2} x_{1} ^k$
$=\sum \limits _{k=m-{\rm ord }(p_{11})+1}^{{\rm ord }(p_{11})-1} g_{0,k}^{m} x_{1} ^{m- k}x_{2} x_{1} ^k$.

Thus equation  holds when $l=0$. Assume ${\rm ord }(p_{11})-m_{12}-1>l>0$. See

$\overline{l}_{1}^{m}[2]^{-} =
\sum \limits _{k=m-{\rm ord }(p_{11})+l+1}^{{\rm ord }(p_{11})-1} g_{l,k}^{m} x_{1} ^{m- k}x_{2} x_{1} ^k
$

$=
\sum \limits _{k=m-{\rm ord }(p_{11})+l+2}^{{\rm ord }(p_{11})-1} g_{l,k}^{m} x_{1} ^{m- k}x_{2} x_{1} ^k
+ g_{l,m-{\rm ord }(p_{11})+l+1}^{m} x_{1} ^{{\rm ord }(p_{11})-l-1}x_{2} x_{1} ^{m-{\rm ord }(p_{11})+l+1}$

\hfill (${\rm ord }(p_{11})-m_{12}-1>l$ and {\rm (ii)})

$=
\sum \limits _{k=m-{\rm ord }(p_{11})+l+2}^{{\rm ord }(p_{11})-1} g_{l,k}^{m} x_{1} ^{m- k}x_{2} x_{1} ^k
$

\hfill $+ g_{l,m-{\rm ord }(p_{11})+l+1}^{m}
(-\sum \limits _{k=1}^{{\rm ord }(p_{11})-l-1}f_{{\rm ord }(p_{11})-l-1,p_{12}}^{k}x_{1} ^{{\rm ord }(p_{11})-l-1-k}x_{2} x_{1} ^{k}) x_{1} ^{m-{\rm ord }(p_{11})+l+1}$

{\scriptsize$=
\sum \limits _{k=m-{\rm ord }(p_{11})+l+2}^{{\rm ord }(p_{11})-1} g_{l,k}^{m} x_{1} ^{m- k}x_{2} x_{1} ^k
-g_{l,m-{\rm ord }(p_{11})+l+1}^{m}
\sum \limits _{k=1}^{2{\rm ord }(p_{11})-2-m-l}f_{{\rm ord }(p_{11})-l-1,p_{12}}^{k}x_{1} ^{{\rm ord }(p_{11})-l-1-k}x_{2} x_{1} ^{k+m-{\rm ord }(p_{11})+l+1}$}

\hfill ($k+m-{\rm ord }(p_{11})+l+1\leq {\rm ord }(p_{11})-1$)

$=
\sum \limits _{k=m-{\rm ord }(p_{11})+l+2}^{{\rm ord }(p_{11})-1} g_{l,k}^{m} x_{1} ^{m- k}x_{2} x_{1} ^k
-g_{l,m-{\rm ord }(p_{11})+l+1}^{m}
\sum \limits _{k=m-{\rm ord }(p_{11})+l+2}^{{\rm ord }(p_{11})-1}f_{{\rm ord }(p_{11})-l-1
,p_{12}}
^{k-m+{\rm ord }(p_{11})-l-1}x_{1} ^{m- k}x_{2} x_{1} ^{k}$

$=
\sum \limits _{k=m-{\rm ord }(p_{11})+l+2}^{{\rm ord }(p_{11})-1} \{g_{l,k}^{m}
-g_{l,m-{\rm ord }(p_{11})+l+1}^{m}f_{{\rm ord }(p_{11})-l-1,p_{12}}
^{k-m+{\rm ord }(p_{11})-l-1}\} x_{1} ^{m- k}x_{2} x_{1} ^k
$

$=
\sum \limits _{k=m-{\rm ord }(p_{11})+l+2}^{{\rm ord }(p_{11})-1} g_{l+1,k}^{m} x_{1} ^{m- k}x_{2} x_{1} ^k$.
Consequently, equation holds.

{\rm (iv)} It is clear by  {\rm (iii)}.
\hfill $\Box$

\begin{Lemma}\label{2.5'} For $\forall\ 0\leq t<m_{12}$.

{\rm (i)} $\sum \limits _{i=0}^{m_{12}-t}f_{m_{12}-t,p_{11}^{-m_{12}+t+1}}^{i}=0$,
$\sum \limits _{i=0}^{m_{21}-t}\bar{f}_{m_{21}-t,p_{22}^{-m_{21}+t+1}}^{i}=0$.

{\rm (ii)} $\sum \limits _{i=0}^{p}f_{m_{12}-t,p_{11}^{-m_{12}+t+1}}^{i}=p_{11}^{-p}f_{m_{12}-t-1,
p_{11}^{-m_{12}+t+2}}^{p}$
for $\forall\ 0\leq p<m_{12}-t$. $\sum \limits _{i=0}^{p}\bar{f}_{m_{21}-t,p_{22}^{-m_{21}+t+1}}^{i}=p_{22}^{-p}\bar{f}_{m_{21}-t-1,
p_{22}^{-m_{21}+t+2}}^{p}$
for $\forall\ 0\leq p<m_{21}-t$.
\end {Lemma}
\noindent {\it Proof.}
{\rm (i)} $\sum \limits _{i=0}^{m_{12}-t}f_{m_{12}-t,p_{11}^{-m_{12}+t+1}}^{i}$
$=\prod \limits _{i=0}^{m_{12}-t-1}(1-p_{11}^{i}p_{11}^{-m_{12}+t+1})=0$
by Lemma \ref {1.2}{\rm (i)} .

{\rm (ii)} We know $\sum \limits _{i=0}^{m_{12}-t-1}f_{m_{12}-t,p_{11}^{-m_{12}+t+1}}^{i}
=-f_{m_{12}-t,p_{11}^{-m_{12}+t+1}}^{m_{12}-t}
$
$=-(-1)^{m_{12}-t}p_{11}^{\frac{(m_{12}-t)(m_{12}-t-1)}{2}} $

\noindent $(p_{11}^{-m_{12}+t+1})^{m_{12}-t}
=p_{11}^{-(m_{12}-t-1)}(-1)^{m_{12}-t-1}p_{11}^{\frac{(m_{12}-t-1)(m_{12}-t-2)}{2}} (p_{11}^{-m_{12}+t+2})^{m_{12}-t-1}=p_{11}^{-(m_{12}-t-1)}
$

\noindent $f_{m_{12}-t-1,p_{11}^{-m_{12}+t+2}}^{m_{12}-t-1}.$

$\sum \limits _{i=0}^{p-1}f_{m_{12}-t,p_{11}^{-m_{12}+t+1}}^{i}=p_{11}^{-p}f_{m_{12}-t-1,p_{11}
^{-m_{12}+t+2}}^{p}-f_{m_{12}-t,p_{11}^{-m_{12}+t+1}}^{p}
$

$=p_{11}^{-p}(-1)^{p}p_{11}^{\frac{p(p-1)}{2}} (p_{11}^{-m_{12}+t+2})^{p}{\scriptsize\left(\begin{array}{cc} m_{12}-t-1\\
p \end{array}\right)_{p_{11}}}-(-1)^{p}p_{11}^{\frac{p(p-1)}{2}} (p_{11}^{-m_{12}+t+1})^{p}{\scriptsize\left(\begin{array}{cc} m_{12}-t\\
p \end{array}\right)_{p_{11}}}$

$=(-1)^{p}p_{11}^{\frac{p(p-1)}{2}} (p_{11}^{-m_{12}+t+1})^{p}({\scriptsize\left(\begin{array}{cc} m_{12}-t-1\\
p \end{array}\right)_{p_{11}}}-{\scriptsize\left(\begin{array}{cc} m_{12}-t\\
p \end{array}\right)_{p_{11}}})$

$=(-1)^{p}p_{11}^{\frac{p(p-1)}{2}} (p_{11}^{-m_{12}+t+1})^{p}(-p_{11}^{m_{12}-t-p}{\scriptsize\left(\begin{array}{cc} m_{12}-t-1\\
p-1 \end{array}\right)_{p_{11}}})$

$=(-1)^{p-1}p_{11}^{\frac{(p-1)(p-2)}{2}}p_{11}^{p-1} (p_{11}^{-m_{12}+t+2})^{(p-1)}p_{11}^{-m_{12}+t+2}p_{11}^{-p}p_{11}^{m_{12}-t}p_{11}^{-p}
{\scriptsize\left(\begin{array}{cc} m_{12}-t-1\\
p-1 \end{array}\right)_{p_{11}}}$

$=p_{11}^{-p+1}f_{m_{12}-t-1,p_{11}
^{-m_{12}+t+2}}^{p-1}$.   \hfill $\Box$

Set $(p_{11}^{m}p_{12}-1)^{{\rm ord }(p_{11})+l}=\alpha_{m,l}^{0}$ for $\forall\ 0\leq m\leq m_{12}$,
$\alpha_{m,l}^{i}-p_{11}^{-i}\alpha_{m+1,l}^{i}
=\alpha_{m,l}^{i+1}$ for $\forall\ 0\leq l,m< m_{12},0\leq i\leq l$, $(p_{22}^{m}p_{21}-1)^{{\rm ord }(p_{22})+l}=\bar{\alpha}_{m,l}^{0}$ for $\forall\ 0\leq m\leq m_{21}$,
$\bar{\alpha}_{m,l}^{i}-p_{22}^{-i}\bar{\alpha}_{m+1,l}^{i}
=\bar{\alpha}_{m,l}^{i+1}$ for $\forall\ 0\leq l,m< m_{21},0\leq i\leq l$.
It is clear $p_{11}^{m}p_{12}-1\neq p_{11}^{n}p_{12}-1$ for $0\leq m\neq n\leq m_{12}$, $p_{22}^{m}p_{21}-1\neq p_{22}^{n}p_{21}-1$ for $0\leq m\neq n\leq m_{21}$.

\begin {Proposition} \label {2.6}
{\rm (i)} For $0\leq l<m_{12}$,
then
 \begin {eqnarray}  \overline{l}_{1}^{{\rm ord }(p_{11})+l}[2]^{-}=0 &\Longleftrightarrow& \alpha_{k,l}^{l+1}=0
 \hbox { for } \forall\ k=0,1,\ldots,m_{12}-l-1. \label {e1.1}
\end {eqnarray}

{\rm (ii)} For $0\leq l<m_{21}$,
then
 \begin {eqnarray}  \overline{l}_{2}^{{\rm ord }(p_{22})+l}[1]^{-}=0 &\Longleftrightarrow& \bar{\alpha}_{k,l}^{l+1}=0
 \hbox { for } \forall\ k=0,1,\ldots,m_{21}-l-1. \label {e1.1''}
\end {eqnarray}
\end {Proposition}
\noindent {\it Proof.} {\rm (i)} We obtain $ \overline{l}_{1}^{{\rm ord }(p_{11})+l}[2]^{-}=0$ if
and only if $g_{{\rm ord }(p_{11})-1-m_{12},{\rm ord }(p_{11})+l-m_{12}+k}^{{\rm ord }(p_{11})+l}=0$ for $\forall\ k=0,1,\ldots,m_{12}-l-1$.

$g_{{\rm ord }(p_{11})-1-m_{12},{\rm ord }(p_{11})+l-m_{12}+k}^{{\rm ord }(p_{11})+l}
$

$=(-1)^{{\rm ord }(p_{11})+l-m_{12}+k} {\scriptsize\left(\begin{array}{cc} {\rm ord }(p_{11})+l\\
{\rm ord }(p_{11})+l-m_{12}+k \end{array}\right)}
-\sum \limits _{j=0}^{{\rm ord }(p_{11})-2-m_{12}}(-1)^{l+1+j} {\scriptsize\left(\begin{array}{cc} {\rm ord }(p_{11})+l\\
l+1+j \end{array}\right)}$

$(-1)^{k-1}p_{11}^{\frac{(k+1)k}{2}}p_{12}^{{\rm ord }(p_{11})-1-m_{12}+k-j}
{\scriptsize\left(\begin{array}{cc} {\rm ord }(p_{11})-1-j\\
{\rm ord }(p_{11})-1-m_{12}+k-j \end{array}\right)_{p_{11}}}{\scriptsize\left(\begin{array}{cc} {\rm ord }(p_{11})-1-m_{12}+k-j-1\\
{\rm ord }(p_{11})-1-m_{12}-j-1 \end{array}\right)_{p_{11}}}$

$
=(-1)^{{\rm ord }(p_{11})+l-m_{12}+k} {\scriptsize\left(\begin{array}{cc} {\rm ord }(p_{11})+l\\
{\rm ord }(p_{11})+l-m_{12}+k \end{array}\right)}
-\sum \limits _{j=l+1}^{{\rm ord }(p_{11})-m_{12}+l-1}(-1)^{j} {\scriptsize\left(\begin{array}{cc} {\rm ord }(p_{11})+l\\
j \end{array}\right)}$

$(-1)^{k-1}p_{11}^{\frac{(k+1)k}{2}}p_{12}^{{\rm ord }(p_{11})-m_{12}+k+l-j}
{\scriptsize\left(\begin{array}{cc} {\rm ord }(p_{11})+l-j\\
{\rm ord }(p_{11})-m_{12}+k+l-j \end{array}\right)_{p_{11}}}{\scriptsize\left(\begin{array}{cc} {\rm ord }(p_{11})-m_{12}+k+l-j-1\\
{\rm ord }(p_{11})-m_{12}+l-j-1 \end{array}\right)_{p_{11}}}$

$
=(-1)^{{\rm ord }(p_{11})+l-m_{12}+k} {\scriptsize\left(\begin{array}{cc} {\rm ord }(p_{11})+l\\
{\rm ord }(p_{11})+l-m_{12}+k \end{array}\right)}
-\sum \limits _{j=l+1}^{{\rm ord }(p_{11})-m_{12}+l-1}(-1)^{j} {\scriptsize\left(\begin{array}{cc} {\rm ord }(p_{11})+l\\
j \end{array}\right)}$

\hfill $(-1)^{k-1}p_{11}^{\frac{(k+1)k}{2}}p_{12}^{{\rm ord }(p_{11})-m_{12}+k+l-j}
\frac{({\rm ord }(p_{11})+l-j)_{p_{11}}({\rm ord }(p_{11})+l-j-1)_{p_{11}}\cdots ({\rm ord }(p_{11})+l-j-m_{12})_{p_{11}}}{(m_{12}-k)_{p_{11}}^{!}(k)_{p_{11}}^{!}(-m_{12}+k+{\rm ord }(p_{11})+l-j)_{p_{11}}}$

$
=(-1)^{{\rm ord }(p_{11})+l-m_{12}+k} {\scriptsize\left(\begin{array}{cc} {\rm ord }(p_{11})+l\\
{\rm ord }(p_{11})+l-m_{12}+k \end{array}\right)}
-\sum \limits _{j=l+1}^{{\rm ord }(p_{11})-m_{12}+l-1}(-1)^{j} {\scriptsize\left(\begin{array}{cc} {\rm ord }(p_{11})+l\\
j \end{array}\right)}(-1)^{k-1}p_{11}^{\frac{(k+1)k}{2}}$

$p_{12}^{{\rm ord }(p_{11})-m_{12}+k+l-j}\frac{1}{(m_{12}-k)_{p_{11}}^{!}(k)_{p_{11}}^{!}(1-p_{11})^{m_{12}}}
\frac{(1-p_{11}^{{\rm ord }(p_{11})+l-j})(1-p_{11}^{{\rm ord }(p_{11})+l-j-1})\cdots (1-p_{11}^{{\rm ord }(p_{11})+l-j-m_{12}})}{(1-p_{11}^{-m_{12}+k+{\rm ord }(p_{11})+l-j})}$

$
=(-1)^{{\rm ord }(p_{11})+l-m_{12}+k} {\scriptsize\left(\begin{array}{cc} {\rm ord }(p_{11})+l\\
{\rm ord }(p_{11})+l-m_{12}+k \end{array}\right)}
+\sum \limits _{j=l+1}^{{\rm ord }(p_{11})-m_{12}+l-1}(-1)^{j} {\scriptsize\left(\begin{array}{cc} {\rm ord }(p_{11})+l\\
j \end{array}\right)}(-1)^{k}p_{11}^{\frac{(k+1)k}{2}}$

$p_{12}^{{\rm ord }(p_{11})-m_{12}+k+l-j}\frac{1}{(m_{12}-k)_{p_{11}}^{!}(k)_{p_{11}}^{!}(1-p_{11})^{m_{12}}}
(1-p_{11}^{{\rm ord }(p_{11})+l-j})\cdots (1-p_{11}^{-m_{12}+k+{\rm ord }(p_{11})+l-j+1})$

\hfill $(1-p_{11}^{-m_{12}+k+{\rm ord }(p_{11})+l-j-1})\cdots(1-p_{11}^{{\rm ord }(p_{11})+l-j-m_{12}})$

$
=\sum \limits _{j=0}^{{\rm ord }(p_{11})+l}(-1)^{j} {\scriptsize\left(\begin{array}{cc} {\rm ord }(p_{11})+l\\
j \end{array}\right)}(-1)^{k}p_{11}^{\frac{(k+1)k}{2}}p_{12}^{{\rm ord }(p_{11})-m_{12}+k+l-j}\frac{1}{(m_{12}-k)_{p_{11}}^{!}(k)_{p_{11}}^{!}(1-p_{11})^{m_{12}}}$

$
(1-p_{11}^{{\rm ord }(p_{11})+l-j})\cdots (1-p_{11}^{-m_{12}+k+{\rm ord }(p_{11})+l-j+1})(1-p_{11}^{-m_{12}+k+{\rm ord }(p_{11})+l-j-1})\cdots(1-p_{11}^{{\rm ord }(p_{11})+l-j-m_{12}})$

$
=\sum \limits _{j=0}^{{\rm ord }(p_{11})+l}(-1)^{j} {\scriptsize\left(\begin{array}{cc} {\rm ord }(p_{11})+l\\
j \end{array}\right)}(-1)^{k}p_{11}^{\frac{(k+1)k}{2}}p_{12}^{{\rm ord }(p_{11})-m_{12}+k+l-j}\frac{1}{(m_{12}-k)_{p_{11}}^{!}(k)_{p_{11}}^{!}(1-p_{11})^{m_{12}}}$

\hfill $\prod \limits _{t=0}^{m_{12}-k-1}(1-p_{11}^{t}p_{11}^{-m_{12}+k+{\rm ord }(p_{11})+l-j+1})
\prod \limits _{t=0}^{k-1}(1-p_{11}^{t}p_{11}^{{\rm ord }(p_{11})+l-j-m_{12}})$

$
=\sum \limits _{j=0}^{{\rm ord }(p_{11})+l}(-1)^{j} {\scriptsize\left(\begin{array}{cc} {\rm ord }(p_{11})+l\\
j \end{array}\right)}(-1)^{k}p_{11}^{\frac{(k+1)k}{2}}p_{12}^{{\rm ord }(p_{11})-m_{12}+k+l-j}\frac{1}{(m_{12}-k)_{p_{11}}^{!}(k)_{p_{11}}^{!}(1-p_{11})^{m_{12}}}$

\hfill $\sum \limits _{t=0}^{m_{12}-k}f_{m_{12}-k,p_{11}^{-m_{12}+k+{\rm ord }(p_{11})+l-j+1}}^{t}
\sum \limits _{t=0}^{k}f_{k,p_{11}^{{\rm ord }(p_{11})+l-j-m_{12}}}^{t}$

$
=\sum \limits _{j=0}^{{\rm ord }(p_{11})+l}(-1)^{j} {\scriptsize\left(\begin{array}{cc} {\rm ord }(p_{11})+l\\
j \end{array}\right)}(-1)^{k}p_{11}^{\frac{(k+1)k}{2}}p_{12}^{{\rm ord }(p_{11})-m_{12}+k+l-j}\frac{1}{(m_{12}-k)_{p_{11}}^{!}(k)_{p_{11}}^{!}(1-p_{11})^{m_{12}}}$

\hfill $\sum \limits _{t_{1}=0}^{m_{12}-k}(-1)^{t_{1}}p_{11}^{\frac{t_{1}(t_{1}-1)}{2}} (p_{11}^{-m_{12}+k+{\rm ord }(p_{11})+l-j+1})^{t_{1}}{\scriptsize\left(\begin{array}{cc} m_{12}-k\\
t_{1} \end{array}\right)_{p_{11}}}$

\hfill $
\sum \limits _{t_{2}=0}^{k}(-1)^{t_{2}}p_{11}^{\frac{t_{2}(t_{2}-1)}{2}} (p_{11}^{{\rm ord }(p_{11})+l-j-m_{12}})^{t_{2}}{\scriptsize\left(\begin{array}{cc} k\\
t_{2} \end{array}\right)_{p_{11}}}$

$
=(-1)^{k}p_{11}^{\frac{(k+1)k}{2}}\frac{1}{(m_{12}-k)_{p_{11}}^{!}(k)_{p_{11}}^{!}(1-p_{11})^{m_{12}}}
p_{12}^{-m_{12}+k}$

\hfill $\sum \limits _{t_{1}=0}^{m_{12}-k}(-1)^{t_{1}}p_{11}^{\frac{t_{1}(t_{1}-1)}{2}}
(p_{11}^{-m_{12}+k+1})^{t_{1}}{\scriptsize\left(\begin{array}{cc} m_{12}-k\\
t_{1} \end{array}\right)_{p_{11}}}
\sum \limits _{t_{2}=0}^{k}(-1)^{t_{2}}p_{11}^{\frac{t_{2}(t_{2}-1)}{2}} (p_{11}^{-m_{12}})^{t_{2}}{\scriptsize\left(\begin{array}{cc} k\\
t_{2} \end{array}\right)_{p_{11}}}
$

\hfill $\sum \limits _{j=0}^{{\rm ord }(p_{11})+l}(-1)^{j} {\scriptsize\left(\begin{array}{cc} {\rm ord }(p_{11})+l\\
j \end{array}\right)}p_{12}^{{\rm ord }(p_{11})+l-j}(p_{11}^{{\rm ord }(p_{11})+l-j})^{t_{1}}(p_{11}^{{\rm ord }(p_{11})+l-j})^{t_{2}}$

$
=(-1)^{k}p_{11}^{\frac{(k+1)k}{2}}\frac{1}{(m_{12}-k)_{p_{11}}^{!}
(k)_{p_{11}}^{!}(1-p_{11})^{m_{12}}}
p_{12}^{-m_{12}+k}\sum \limits _{t_{1}=0}^{m_{12}-k}f_{m_{12}-k,p_{11}
^{-m_{12}+k+1}}^{t_{1}}
\sum \limits _{t_{2}=0}^{k}f_{k,p_{11}
^{-m_{12}}}^{t_{2}}\alpha_{t_{1}+t_{2},l}^{0}
$.

We obtain $g_{{\rm ord }(p_{11})-1-m_{12},{\rm ord }(p_{11})+l-m_{12}+k}^{{\rm ord }(p_{11})+l}=0$
if and only if $\sum \limits _{t_{1}=0}^{m_{12}-k}f_{m_{12}-k,p_{11}
^{-m_{12}+k+1}}^{t_{1}}
$

\noindent $\sum \limits _{t_{2}=0}^{k}f_{k,p_{11}
^{-m_{12}}}^{t_{2}}\alpha_{t_{1}+t_{2},l}^{0}=0$ for
$\forall\ k=0,1,\ldots,m_{12}-l-1$.

We know $\sum \limits _{t_{1}=0}^{m_{12}-k}f_{m_{12}-k,p_{11}
^{-m_{12}+k+1}}^{t_{1}}\sum \limits _{t_{2}=0}^{k}f_{k,p_{11}
^{-m_{12}}}^{t_{2}}\alpha_{t_{1}+t_{2},l}^{0}
$

\noindent $=\sum \limits _{t_{1}=0}^{m_{12}-k-p}
\sum \limits _{t_{2}=0}^{k}a^{-pt_{1}}f_{m_{12}-k-p,p_{11}
^{-m_{12}+k+p+1}}^{t_{1}}
f_{k,p_{11}
^{-m_{12}}}^{t_{2}}\alpha_{t_{1}+t_{2},l}^{p}$
for $1\leq p\leq m_{12}-k$:

$\sum \limits _{t_{1}=0}^{m_{12}-k}
\sum \limits _{t_{2}=0}^{k}f_{m_{12}-k,p_{11}
^{-m_{12}+k+1}}^{t_{1}}f_{k,p_{11}
^{-m_{12}}}^{t_{2}}\alpha_{t_{1}+t_{2},l}^{0}
$

$=\sum \limits _{t_{1}=0}^{m_{12}-k-1}
\sum \limits _{t_{2}=0}^{k}\sum \limits _{q=0}^{t_{1}}f_{m_{12}-k,p_{11}
^{-m_{12}+k+1}}^{q}f_{k,p_{11}
^{-m_{12}}}^{t_{2}}
(\alpha_{t_{1}+t_{2},l}^{0}-\alpha_{t_{1}+t_{2}+1,l}^{0})
$

$+
\sum \limits _{t_{2}=0}^{k}\sum \limits _{q=0}^{m_{12}-k}f_{m_{12}-k,p_{11}
^{-m_{12}+k+1}}^{q}f_{k,p_{11}
^{-m_{12}}}^{t_{2}}
\alpha_{m_{12}-k+t_{2},l}^{0}$

$=\sum \limits _{t_{1}=0}^{m_{12}-k-1}
\sum \limits _{t_{2}=0}^{k}a^{-t_{1}}f_{m_{12}-k-1,p_{11}
^{-m_{12}+k+2}}^{t_{1}}f_{k,p_{11}
^{-m_{12}}}^{t_{2}}
\alpha_{t_{1}+t_{2},l}^{1}
$ by Lemma \ref {2.5'}. It is hold for $p=1$.

$\sum \limits _{t_{1}=0}^{m_{12}-k-p}
\sum \limits _{t_{2}=0}^{k}a^{-pt_{1}}f_{m_{12}-k-p,p_{11}
^{-m_{12}+k+p+1}}^{t_{1}}
f_{k,p_{11}
^{-m_{12}}}^{t_{2}}\alpha_{t_{1}+t_{2},l}^{p}$

$=\sum \limits _{t_{1}=0}^{m_{12}-k-p-1}
\sum \limits _{t_{2}=0}^{k}a^{-pt_{1}}\sum \limits _{q=0}^{t_{1}}f_{m_{12}-k-p,p_{11}
^{-m_{12}+k+p+1}}^{q}
f_{k,p_{11}
^{-m_{12}}}^{t_{2}}(\alpha_{t_{1}+t_{2},l}^{p}-a^{-p}\alpha_{t_{1}+t_{2}+1,l}^{p})
$

$+
\sum \limits _{t_{2}=0}^{k}a^{-p(m_{12}-k-p)}\sum \limits _{q=0}^{m_{12}-k-p}f_{m_{12}-k-p,p_{11}
^{-m_{12}+k+p+1}}^{q}
f_{k,p_{11}
^{-m_{12}}}^{t_{2}}\alpha_{m_{12}-k-p+t_{2},l}^{p}$

$=\sum \limits _{t_{1}=0}^{m_{12}-k-p-1}
\sum \limits _{t_{2}=0}^{k}a^{-(p+1)t_{1}}f_{m_{12}-k-p-1,p_{11}
^{-m_{12}+k+p+2}}^{t_{1}}
f_{k,p_{11}
^{-m_{12}}}^{t_{2}}\alpha_{t_{1}+t_{2},l}^{p+1}
$.

Then we obtain
$\sum \limits _{t_{2}=0}^{k}
f_{k,p_{11}
^{-m_{12}}}^{t_{2}}\alpha_{t_{2},l}^{m_{12}-k}=0$ for
$\forall\ k=0,1,\ldots,m_{12}-l-1$.

We know $\alpha_{t_{2},l}^{m_{12}-k}=0$ for $\forall\ 0\leq t_{2}\leq k,\forall\ k=0,1,\ldots,m_{12}-l-1$:

If $k=0$, $\sum \limits _{t_{2}=0}^{k}
f_{k,p_{11}
^{-m_{12}}}^{t_{2}}\alpha_{t_{2},l}^{m_{12}-k}=\alpha_{0,l}^{m_{12}-0}=0$, it is clear.

Assumed that $\leq k$ is hold, now consider $k+1$.
$\sum \limits _{t_{2}=0}^{k+1}
f_{k+1,p_{11}
^{-m_{12}}}^{t_{2}}\alpha_{t_{2},l}^{m_{12}-k-1}=0$, we obtain
$\alpha_{t_{2},l}^{m_{12}-k-1}=p_{11}^{-(m_{12}-k-1)}\alpha_{t_{2}+1,l}^{m_{12}-k-1}$
by $\alpha_{t_{2},l}^{m_{12}-k}=0$.
Then $\alpha_{t_{2},l}^{m_{12}-k-1}=p_{11}^{-(k+1-t_{2})(m_{12}-k-1)}\alpha_{k+1,l}^{m_{12}-k-1}$
 for $\forall\ 0\leq t_{2}\leq k$. Thus $\sum \limits _{t_{2}=0}^{k+1}
f_{k+1,p_{11}
^{-m_{12}}}^{t_{2}}\alpha_{t_{2},l}^{m_{12}-k-1}=\sum \limits _{t_{2}=0}^{k+1}
f_{k+1,p_{11}
^{-m_{12}}}^{t_{2}}p_{11}^{-(k+1-t_{2})(m_{12}-k-1)}\alpha_{k+1,l}^{m_{12}-k-1}=0$.
It is clear $\sum \limits _{t_{2}=0}^{k+1}
f_{k+1,p_{11}
^{-m_{12}}}^{t_{2}}p_{11}^{-(k+1-t_{2})(m_{12}-k-1)}$

$=\sum \limits _{t_{2}=0}^{k+1}
(-1)^{t_{2}}p_{11}^{\frac{t_{2}(t_{2}-1)}{2}} (p_{11}^{-m_{12}})^{t_{2}}{\scriptsize\left(\begin{array}{cc} k+1\\
t_{2} \end{array}\right)_{p_{11}}}p_{11}^{-(k+1-t_{2})(m_{12}-k-1)}$

$=\sum \limits _{t_{2}=0}^{k+1}
(-1)^{t_{2}}p_{11}^{\frac{t_{2}(t_{2}-1)}{2}} p_{11}^{-(k+1)(m_{12}-k-1)}
p_{11}^{t_{2}(-k-1)}
{\scriptsize\left(\begin{array}{cc} k+1\\
t_{2} \end{array}\right)_{p_{11}}}$

$=p_{11}^{-(k+1)(m_{12}-k-1)}\prod \limits _{t_{2}=0}^{k}(1-p_{11}^{t_{2}}p_{11}^{(-k-1)})\neq0$
by Lemma \ref {1.2}{\rm (i)} and $k<m_{12}<{\rm ord }(p_{11})-1$.

Then $\alpha_{k+1,l}^{m_{12}-k-1}=0$, then
$\alpha_{t_{2},l}^{m_{12}-k-1}=0$ for $\forall\ 0\leq t_{2}\leq k+1$.
\hfill $\Box$

\begin {Theorem} \label {2.6''}
{\rm (i)} $\overline{l}_{1}^{m}[2]^ -=0$ if $m\geq m_{12}+{\rm ord }(p_{11})$;
$\overline{l}_{2}^{m}[1]^ -=0$ if $m\geq m_{21}+{\rm ord }(p_{22})$.

{\rm (ii)} $\overline{l}_{1}^{m}[2]^ -\neq0$ if $m<{\rm ord }(p_{11})$; $\overline{l}_{2}^{m}[1]^ -\neq0$ if $m<{\rm ord }(p_{22})$.

{\rm (iii)} For $0\leq l<m_{12}$, then
 $\overline{l}_{1}^{{\rm ord }(p_{11})+l}[2]^{-}=0$
 $\Longleftrightarrow$
 $ \alpha_{k,l}^{l+1}=0$
 for $\forall\ k=0,1,\ldots,m_{12}-l-1$.
For $0\leq l<m_{21}$,  then
$ \overline{l}_{2}^{{\rm ord }(p_{22})+l}[1]^{-}=0$
$\Longleftrightarrow$
$ \bar{\alpha}_{k,l}^{l+1}=0$
for $\forall\ k=0,1,\ldots,m_{21}-l-1$.
\end {Theorem}
\noindent {\it Proof.} {\rm (i)} We obtain $\overline{l}_{1}^{m_{12}+{\rm ord }(p_{11})}[2]^{-}=\sum \limits _{k=0}^{m_{12}+{\rm ord }(p_{11})} ( - 1)^k {\scriptsize\left(\begin{array}{cc} m_{12}+{\rm ord }(p_{11})\\
k \end{array}\right)} x_{1} ^{m_{12}+{\rm ord }(p_{11})- k}x_{2} x_{1} ^k$
$=0$ by Lemma \ref {6.5}(iii). The other is clear by Lemma \ref {6.5} {\rm (ii)}{\rm (2)}.
\hfill $\Box$

\begin {Example}  \label {2.7}
{\rm (i)} If  $m_{12}=1$, then
$\overline{l}_{1}^{{\rm ord }(p_{11})}[2]^{-}=0 \Longleftrightarrow \alpha_{0,0}^{1}=0$ ,i.e. $(p_{12}-1)^{{\rm ord }(p_{11})}=(p_{11}p_{12}-1)^{{\rm ord }(p_{11})}$.
If  $m_{21}=1$, then
$\overline{l}_{2}^{{\rm ord }(p_{22})}[1]^{-}=0 \Longleftrightarrow \bar{\alpha}_{0,0}^{1}=0$ ,i.e. $(p_{21}-1)^{{\rm ord }(p_{22})}=(p_{22}p_{21}-1)^{{\rm ord }(p_{22})}$.

{\rm (ii)} If  $m_{12}=2$, then
$\overline{l}_{1}^{{\rm ord }(p_{11})}[2]^{-}=0 \Longleftrightarrow \alpha_{0,0}^{1}
=\alpha_{1,0}^{1}=0$ ,i.e. $(p_{12}-1)^{{\rm ord }(p_{11})}=(p_{11}p_{12}-1)^{{\rm ord }(p_{11})}
=(p_{11}^{2}p_{12}-1)^{{\rm ord }(p_{11})}$.
$\overline{l}_{1}^{{\rm ord }(p_{11})+1}[2]^{-}=0 \Longleftrightarrow \alpha_{0,1}^{2}=0$ ,i.e. $(p_{12}-1)^{{\rm ord }(p_{11})+1}-(1+p_{11}^{-1})(p_{11}p_{12}-1)^{{\rm ord }(p_{11})+1}
+p_{11}^{-1}(p_{11}^{2}p_{12}-1)^{{\rm ord }(p_{11})+1}=0$.
If  $m_{21}=2$, then
$\overline{l}_{2}^{{\rm ord }(p_{22})}[1]^{-}=0 \Longleftrightarrow \bar{\alpha}_{0,0}^{1}
=\bar{\alpha}_{1,0}^{1}=0$ ,i.e. $(p_{21}-1)^{{\rm ord }(p_{22})}=(p_{22}p_{21}-1)^{{\rm ord }(p_{22})}
=(p_{22}^{2}p_{21}-1)^{{\rm ord }(p_{22})}$.
$\overline{l}_{2}^{{\rm ord }(p_{22})+1}[1]^{-}=0 \Longleftrightarrow \bar{\alpha}_{0,1}^{2}=0$ ,i.e. $(p_{21}-1)^{{\rm ord }(p_{22})+1}-(1+p_{22}^{-1})(p_{22}p_{21}-1)^{{\rm ord }(p_{22})+1}
+p_{22}^{-1}(p_{22}^{2}p_{21}-1)^{{\rm ord }(p_{22})+1}=0$.
\end {Example}

\begin {Example}  \label {6.10} Assume that $m_{12}={\rm ord }(p_{11})-2$.

{\rm (i)} If $p_{12}=-p_{11}$. Then $\overline{l}_{1}^{m}[2]^-\neq0$ if and only if
$ m_{12}+\ {\rm ord }(p_{11})-2 \geq m\geq0$.

{\rm (ii)} If $p_{12}\neq-p_{11}$. Then $\overline{l}_{1}^{m}[2]^-\neq0$ if and only if
$ m_{12}+\ {\rm ord }(p_{11})-1 \geq m\geq0$.

\end {Example}
\noindent {\it Proof.}
$g_{1,{\rm ord }(p_{11})-1}
^{m_{12}+{\rm ord }(p_{11})-1}=g_{1,{\rm ord }(p_{11})-1}
^{2{\rm ord }(p_{11})-3}=g_{0,{\rm ord }(p_{11})-1}^{2{\rm ord }(p_{11})-3}-g_{0,{\rm ord }(p_{11})-2}^{2{\rm ord }(p_{11})-3}f_{{\rm ord }(p_{11})-1,p_{12}}
^{1}$

{\scriptsize$=( - 1)^{{\rm ord }(p_{11})-1} \left(\begin{array}{cc} 2{\rm ord }(p_{11})-3\\
{\rm ord }(p_{11})-1 \end{array}\right)-( - 1)^{{\rm ord }(p_{11})-2} \left(\begin{array}{cc} 2{\rm ord }(p_{11})-3\\
{\rm ord }(p_{11})-2 \end{array}\right)(-1)p_{12}{\scriptsize\left(\begin{array}{cc} {\rm ord }(p_{11})-1\\
1 \end{array}\right)_{p_{11}}}$}

$=g_{0,{\rm ord }(p_{11})-1}^{2{\rm ord }(p_{11})-3}(1-p_{12}{\scriptsize\left(\begin{array}{cc} {\rm ord }(p_{11})-1\\
1 \end{array}\right)_{p_{11}}})=g_{0,{\rm ord }(p_{11})-1}^{2{\rm ord }(p_{11})-3}(1+p_{12}p_{11}^{-1})
$.

{\rm (i)}
$\overline{l}_{1}^{m_{12}+{\rm ord }(p_{11})-1}[2]^{-}=g_{1,{\rm ord }(p_{11})-1}
^{m_{12}+{\rm ord }(p_{11})-1}x_{1} ^{m_{12}}x_{2} x_{1} ^{{\rm ord }(p_{11})-1}
=0$ by Lemma \ref {6.5}{\rm (i)}.

$g_{1,{\rm ord }(p_{11})-2}
^{m_{12}+{\rm ord }(p_{11})-2}=g_{1,{\rm ord }(p_{11})-2}
^{2{\rm ord }(p_{11})-4}=g_{0,{\rm ord }(p_{11})-2}^{2{\rm ord }(p_{11})-4}-g_{0,{\rm ord }(p_{11})-3}^{2{\rm ord }(p_{11})-4}f_{{\rm ord }(p_{11})-1,p_{12}}
^{1}$

{\scriptsize$=( - 1)^{{\rm ord }(p_{11})-2} \left(\begin{array}{cc} 2{\rm ord }(p_{11})-4\\
{\rm ord }(p_{11})-2 \end{array}\right)-( - 1)^{{\rm ord }(p_{11})-3} \left(\begin{array}{cc} 2{\rm ord }(p_{11})-4\\
{\rm ord }(p_{11})-3 \end{array}\right)(-1)p_{12}{\scriptsize\left(\begin{array}{cc} {\rm ord }(p_{11})-1\\
1 \end{array}\right)_{p_{11}}}$}

$=g_{0,{\rm ord }(p_{11})-2}^{2{\rm ord }(p_{11})-4}\{1+\frac{{\rm ord }(p_{11})-2}{{\rm ord }(p_{11})-3} (-1)p_{12}{\scriptsize\left(\begin{array}{cc} {\rm ord }(p_{11})-1\\
1 \end{array}\right)_{p_{11}}}\}$

$=g_{0,{\rm ord }(p_{11})-2}^{2{\rm ord }(p_{11})-4}\{1+\frac{{\rm ord }(p_{11})-2}{{\rm ord }(p_{11})-3} p_{11}({\rm ord }(p_{11})-1)_{p_{11}}\}$

$=g_{0,{\rm ord }(p_{11})-2}^{2{\rm ord }(p_{11})-4}\{1-\frac{{\rm ord }(p_{11})-2}{{\rm ord }(p_{11})-3}\}$

$=-g_{0,{\rm ord }(p_{11})-2}^{2{\rm ord }(p_{11})-4}\frac{1}{{\rm ord }(p_{11})-3}
\neq0$. Then $\overline{l}_{1}^{m_{12}+{\rm ord }(p_{11})-2}[2]^{-}=g_{1,{\rm ord }(p_{11})-1}
^{m_{12}+{\rm ord }(p_{11})-2}x_{1} ^{m_{12}-1}x_{2} x_{1} ^{{\rm ord }(p_{11})-1}$
$+g_{1,{\rm ord }(p_{11})-2}
^{m_{12}+{\rm ord }(p_{11})-2}x_{1} ^{m_{12}}x_{2} x_{1} ^{{\rm ord }(p_{11})-2}
\neq0$.

{\rm (ii)} is clear by Theorem \ref {2.6''}
{\rm (i)}.
\hfill $\Box$

\section*{Acknowledgement}

Y.Z.Z. was partly supported by the National Natural Science Foundation of China (Grant No. 11775177) and the Australian Research Council through Discovery Project grant DP190101529.

\end {document}